                      \def\version{2 February, 2007}

\documentclass[reqno,twoside,11pt]{amsart}
\usepackage{amsmath}
\usepackage{amssymb}

\newcommand{\floor}[1]{\left\lfloor #1 \right\rfloor}
\newcommand{\ceil}[1]{\left\lceil #1 \right\rceil}

\def\d{\delta}

\def\l{\lambda}

\def\L{\Lambda}

\newfam\Bbbfam
\font\tenBbb=msbm10
\font\sevenBbb=msbm7
\font\fiveBbb=msbm5
\textfont\Bbbfam=\tenBbb
\scriptfont\Bbbfam=\sevenBbb
\scriptscriptfont\Bbbfam=\fiveBbb

\newcommand{\R}     {\mathbb{R}}

\newcommand{\N}     {\mathbb{N}}
\renewcommand{\P}   {\mathbb{P}}

\newcommand{\E}     {\mathbb{E}}

\newcommand{\heap}[2]{\genfrac{}{}{0pt}{}{#1}{#2}}

\def\1{{\mathchoice {1\mskip-4mu\mathrm l}      
{1\mskip-4mu\mathrm l}
{1\mskip-4.5mu\mathrm l} {1\mskip-5mu\mathrm l}}}
\newcommand{\ssup}[1] {{\scriptscriptstyle{({#1}})}}
\def\comment#1{}

\newenvironment{proofsect}[1]
{\vskip0.1cm\noindent{\bf #1.}\hskip0.5cm}

\newtheorem{theorem}{Theorem}[section]
\newtheorem{lemma}[theorem]{Lemma}
\newtheorem{prop}[theorem] {Proposition}
\newtheorem{cor}[theorem]  {Corollary}
\newtheorem{remark}[theorem]  {Remark}

\newcommand{\s} {\sigma}

\renewcommand{\d}{{\rm d}}

\newcommand{\Sym}{\mathfrak{S}}

\newcommand{\dist}{{\operatorname {dist}}}

\newcommand{\tr}{{\operatorname {Tr}\,}}

\newcommand{\Bcal}   {{\mathcal B }}
\newcommand{\Ccal}   {{\mathcal C }}

\newcommand{\Hcal}   {{\mathcal H }}

\newcommand{\Mcal}   {{\mathcal M }}

\newcommand{\Pcal}   {{\mathcal P }}

\newcommand{\Scal}   {{\mathcal S }}

\newcommand{\Xcal}   {{\mathcal X }}

\begin{document}

\title[Large deviations in cycles]{\Large\bf  Large deviations for empirical path measures in cycles of integer partitions}

\author[Stefan Adams]{}

\date{}

\maketitle

\thispagestyle{empty}
\vspace{0.2cm}

\centerline {\sc By Stefan Adams\footnote{Max-Planck Institute for Mathematics in the Sciences, Inselstra{\ss}e 22-26, D-04103 Leipzig, Germany, {\tt adams@mis.mpg.de}}}
\vspace{0.4cm}

\vspace{0.2cm}

\centerline{\small(\version)}
\vspace{.3cm}
\begin{quote}
{\small {\bf Abstract:}} 
Consider a large system of $N$ Brownian motions in $\R^d$ on some fixed time interval $[0,\beta]$ with symmetrised initial-terminal condition. That is, for any $i$, the terminal location of the $i$-th motion is affixed to the initial point of the $\s(i)$-th motion, where $\s$ is a uniformly distributed random permutation of $1,\dots,N$. We integrate over all initial points confined in boxes $ \L\subset\R^d $ with respect to the Lebesgue measure, and we divide by a normalisation (partition function). 

Such systems play an important role in quantum physics in the description of Boson systems at positive temperature $1/\beta$.

In this article, we describe the large-$N$ behaviour of the empirical path measure (the mean of the Dirac measures in the $N$ paths) when $ \L\uparrow\R^d $ and $ N/|\L|\to\rho $. The rate function is given as a variational formula involving a certain entropy functional and a Fenchel-Legendre transform. The entropy term governs the large-$N$ behaviour of discrete shape measures of integer partitions. Any integer partition determines a conjugacy class of permutations of certain cycle structure. 

Depending on the dimension and the density $ \rho $, there is phase transition behaviour for the empirical path measure. For certain parameters (high density, large time horizon) and dimensions $ d\ge 3 $ the empirical path measure is not supported on all paths $ [0,\infty)\to\R^d $ which contain a bridge path of any finite multiple of the time horizon $ [0,\beta] $. For dimensions $ d=1,2 $, and for small densities and small time horizon $ [0,\beta] $ in dimensions $ d\ge 3$, the empirical path measure is supported on those paths. In the first regime a finite fraction of the motions lives in cycles of infinite length.

We outline that this transition leads to an empirical path measure interpretation of {\it Bose-Einstein condensation\/}, known for systems of Bosons.

\end{quote}


\eject

\setcounter{section}{0}
\section{Introduction and Results}
We study the large-$N$ behaviour of a system of $ N $ Brownian motions $ B^{\ssup{1}},\ldots, B^{\ssup{N}} $, with time horizon $ [0,\beta] $ in $ \R^d $ confined in subsets $ \L_N\subset \R^d $, i.e., the behaviour of the system under the symmetrised measure
\begin{equation}\label{symPdef1}
\P_N^{\ssup{\rm sym}}=Z_N^{\ssup{\rm sym}}(\beta)^{-1}\frac{1}{N!}\sum_{\s\in\Sym_N}\int_{\L_N}\d x_1\cdots\int_{\L_N}\d x_N\bigotimes_{i=1}^N\mu_{x_i,x_{\s(i)}}^{\beta,N}.
\end{equation}
Here $ \Sym_N $ is the set of permutations of $ 1,\ldots,N $, $ \mu_{x,y}^{\beta,N} $ is the Brownian bridge measure on the time interval $ [0,\beta] $ with initial point $ x\in\L_N $ and terminal point $ y\in\L_N $ and confinement to stay in $ \L_N $ (see \eqref{nnBBM} below), and $ Z_N^{\ssup{\rm sym}}(\beta) $ is the normalisation
\begin{equation}\label{partitionfunction}
Z_N^{\ssup{\rm sym}}(\beta)=\frac{1}{N!}\sum_{\s\in\Sym_N}\int_{\L_N}\d x_1\cdots\int_{\L_N}\d x_N\,\bigotimes_{i=1}^N\mu_{x_i,x_{\s(i)}}^{\beta,N}(\Omega_\beta^N),
\end{equation}
where $ \Omega_\beta $ is the set of continuous functions $ [0,\beta]\to\R^d $.
Hence, the terminal location of the $i$-th motion is affixed to the initial location of the $\sigma(i)$-th motion, where $\sigma$ is a uniformly distributed random permutation. That is, in \eqref{symPdef1} we have two random mechanisms. First we pick uniformly a permutation and after that we pick $ N$ initial points in $ \L_N $ which are permuted according to the chosen permutation to obtain $ N $ terminal points. Then these $ N $ initial and terminal points determine the $ N $ random processes. Finally we average over all permutations and integrate over all initial points in the set $ \L_N $.
Beside the fact that the symmetrised measure $ \P^{\ssup{\rm sym}}_{N} $ is itself of interest,
there are two main motivations in studying the symmetrised measure.

The symmetrisation in \eqref{symPdef1} is described by the set of $ N $ pairs $ (x_1,\ldots,x_N;x_{\sigma(1)},\ldots,x_{\sigma(N)})$ for any permutation $\sigma\in\Sym_N $ and any $ x\in\R^{dN} $. The mixing procedure for the second entry in these pairs has been studied both in \cite{DZ92} and \cite{T02}, which were motivated from asymptotic questions about exchangeable vectors of random variables.  \cite{DZ92} studies large deviations for the empirical measures $\frac 1N\sum_{i=1}^N\delta_{Y_i}$, where $Y_1,\dots,Y_N$ have distribution $\int_\Theta \mu(\d \theta)\, P_N^{\ssup \theta}$ for some distribution $\mu$ on some compact space $\Theta$, and the empirical measures are assumed to satisfy a large deviation principle under $P_N^{\ssup \theta}$ for each $\theta$. In \cite{T02}, a similar problem is studied: given a sequence of random vectors $(Y_1^{\ssup N},\dots,Y_N^{\ssup N})$ such that the empirical measures $\frac 1N\sum_{i=1}^N\delta_{Y_i^{\ssup N}}$ satisfy a large deviation principle, another principle is established for the process of empirical measures $\frac 1N\sum_{i=1}^{\lfloor tN\rfloor}\delta_{X_{i}^{\ssup N}}$, where 
$$
\big(X_1^{\ssup N},\dots,X_N^{\ssup N}\big)=\frac1{N!}\sum_{\s\in\Sym_N}\big(Y_{\s(1)}^{\ssup N},\dots,Y_{\s(N)}^{\ssup N}\big).
$$

Our second main motivation for studying the symmetrised measure $ \P_{N}^{{\ssup{\rm sym}}} $ stems from the applications of Feynman-Kac formulae to express thermodynamic functions in quantum statistical mechanics. These thermodynamic functions are given as traces over exponentials of the Hamilton operator describing the quantum system. There exist two kinds of elementary particles in nature, the {\it Fermions\/} and the {\it Bosons\/}. The state of a system of $ N $ {\it Bosons\/} is described by a symmetrisation procedure like in \eqref{symPdef1}, whereas the state for {\it Fermions\/} is given with the corresponding anti-symmetrisation procedure. Thus one is lead to employ large deviation technique to study the large $ N$-limit for expectations with respect to the symmetrised distribution.

Let $ \L_N\subset\R^d $ a sequence of subsets with $ N/|\L_N|\to\rho \in (0,\infty) $ as $ N\to\infty $. We are going to study large deviations for the empirical path measure
\begin{equation}
L_N=\frac{1}{N}\sum_{i=1}^N\delta_{B^{\ssup{i}}},
\end{equation}
which we conceive as a $ \Omega_\beta^N $-measurable random measure in $ \Mcal_1(\Omega) $, where $ \Omega $ is the set of all continuous paths $[0,\infty)\to\R^d $ (see \eqref{delemp} below for details), under the symmetrised measure $ \P_N^{\ssup{\rm sym}} $. More precisely, we derive a large deviations principle for the distributions of $L_N$ under $\P_{N}^{\ssup{{\rm sym}}}$ as $N\to\infty$ (Theorem~\ref{THMLDP-LN}). (In Section~\ref{sec-ldp} below we recall the notion of a large deviations principle.)



To prove large deviations principles under the symmetrised measure \eqref{symPdef1} one cannot  proceed as in the G\"artner-Ellis-Theorem. A first approach is \cite{AK06} and \cite{AD06} for a different symmetrised measure. However, to see the peculiar correlations due to the symmetrisation one has to study the cycle structure for any permutation. The cycle structure allows to concatenate different Brownian motions to Brownian bridges with larger time horizon. Hence,  a large deviations principle for the empirical path measure as  a random element in the set of probability measures on paths with time horizon $ [0,\infty) $ can be obtained. 

The cycle structure allows to replace the sum over permutations by a sum of integer partitions. 
For any integer $ N $, a partition $ \l $ of $ N $ is the collection of integers $ n_1\ge n_2\ge\cdots\ge n_k\ge 1 , k\in\{1,\ldots,N\}$, such that $ \sum_{i=1}^k n_i=N $. We denote the set of all partitions of $ N $ by $ \Pcal_N $.  Any partition $ \l\in\Pcal_N $ is determined by the sequence $ \{r_k\}_{k=1}^N $ of positive integers $ r_k $ such that $ \sum_{k=1}^Nkr_k =N $, where we write $ r_k(\l)=r_k $. We call the number $ r_k $ an {\it occupation number\/} of the partition.  A cycle of length $ k $ is a chain of permutations, such as $ 1 $ goes to $2$, $2$ goes to $ 3 $, $ 3 $ goes to $4$, etc. until $ k-1 $ goes to $ k $ and finally $ k $ goes to $ 1$. A permutation with exactly $ r_k $ cycles of length $ k $ is said to be of type $ \{r_k\}_{k=1}^N $. Hence, each partition $ \l\in\Pcal_N $ corresponds to a conjugacy class of permutations, i.e., those of the same type, with 
$$ 
\frac{N!}{\prod_{k=1}^N r_k! k^{r_k}} 
$$ elements. The conjugacy classes are the ones where the quantum mechanical trace operation is constant, see Section~\ref{qm-sec} for details.

The rate function of the large deviations principle in Theorem~\ref{THMLDP-LN} is given as a variational formula involving an entropy functional and a Fenchel-Legendre transform. The entropy term governs the large-$N$ behaviour of discrete shape measures of integer partitions. 

In Theorem~\ref{THMfree} we show the thermodynamic limit of the normalisation $ Z_N^{\ssup{\rm sym}}(\beta) $, and provide a variational formula for it. We analyse the variational problem, i.e., the minimisers for the entropy functional, in Theorem~\ref{THM-variational}. In particular, we derive a variational expression for the specific free energy. A phase transition, i.e., a singularity of the free energy, exists for dimensions $ d\ge 3 $ depending on the density $ \rho $ and $ \beta $. For high density, respectively for long time horizons $ [0,\beta] $, the specific free energy is independent of the density. This is the {\it Bose-Einstein condensation\/} for a non-interacting system of Bosons in the thermodynamic limit (see Appendix~\ref{appendix-Bosefunction}).

Our main results, the large deviations principle for $ L_N $ under $ \P_N^{\ssup{\rm sym}} $, in Theorem~\ref{THMLDP-LN} and the analysis of the rate function in Theorem~\ref{THManalysisRF}, give a phase transition for the empirical path measure. 
Let the set 
$$
A_k=\{\omega\otimes_{k\beta}\xi\colon\omega\in\Omega_k,\omega(0)=\omega(k\beta),\xi\in\Omega\},
$$ where $ \Omega_k $ is the set of paths $ [0,k\beta]\to\R^d $ and $ \otimes_{k\beta} $ is the splice of two paths, be given. The empirical path measure for dimensions $ d=1,2, $ and $ \rho<\infty $, or $ \rho<\rho_c $ for $ d\ge 3 $  has support on paths in any $ A_k $, where one can insert any finite number of Brownian motions with time horizon $ [0,\beta] $, i.e., for any $ k\in\N $ one can concatenate exactly $ k $ Brownian motions to paths $ \omega\in A_k $. This is due to the cycle structure of the permutations and the Lebesgue integration over all initial points in the definition of the symmetrised measure $ \P_N^{\ssup{\rm sym}} $. If the density $ \rho $ is high enough for $ d\ge 3 $, i.e., $ \rho>\rho_c $ (or equivalently, if the time horizon is sufficiently large for given density, i.e., $ \beta>\beta_c $, for $ d\ge 3 $), the mean path measure has positive weight on paths with an infinite time horizon, that is, concatenation of  any finite number of Brownian motions with time horizon $ [0,\beta] $ to obtain paths in the $ A_k $'s, is not sufficient.  Then an excess density $ (\rho-\rho_c) $ of Brownian motions with time horizon $ [0,\beta] $ concatenate to infinite long cycles.  


Hence, we have an empirical path measure interpretation of Bose-Einstein condensation.
This novel interpretation adds to the existing notions of Bose-Einstein condensation (\cite{FW71},\cite{S02},\cite{Uel06}), and allows to study systems of interacting Brownian motions. Future work will be devoted to this case \cite{A07a}. Interacting Brownian motions in trap potentials have been so far analysed without symmetrisation, in particular, finite systems for vanishing temperature in \cite{ABK04} and large systems of interacting motions for fixed positive temperature in \cite{ABK05}.

Let us make further remarks on related literature. 

An important work combining combinatorics and large deviations for symmetrised measures is \cite{Toth90}. T\'oth \cite{Toth90} considers $ N $ continuous-time simple random walks on a complete graph with $ \rho N $ vertices, where $ \rho\in (0,1) $ is fixed. He looks at the symmetrised distribution as in \eqref{symPdef1} and adds an exclusion constraint: there is no collision of any two particles during the time interval $ [0,\beta] $. The combinatorial structure of this model enabled him to express the free energy in terms of a cleverly chosen Markov process on $ \N_0 $. Using Freidlin-Wentzell theory, he derives an explicit formula for the large-$ N $ asymptotic of the free energy; in particular he obtains a phase-transition, called {\it Bose-Einstein-condensation\/}, for large $ \beta $ and sufficiently large $ \rho $.

Scaling limits for shape measures of integer partitions in $ \Pcal_N $ under uniform distribution 
are obtained in \cite{Ver96}. Large deviations from this limit behaviour are in \cite{DVZ00}, where large deviations principles for scaled shape measures for partitions as well for strict partitions under uniform distributions are derived. Motivated by the statistics of combinatorial partitions, illustrated by Vershik  in \cite{Ver96}, Benfatto {\it et al.} derived in \cite{BCMP05} limit theorems for statistics of combinatorial partitions for the case of a mean field Bose gas in the grandcanonical ensemble. Here, in contrast to the canonical ensemble, only the mean of the particle number is fixed. \cite{BCMP05} uses Fourier analysis of the corresponding traces to derive a complete description of the statistics of short and long cycles.
For a perturbed mean-field model the density of long cycles for a perturbed mean-field model is analysed in \cite{DMP05}.

In Section~\ref{sec-Results} we describe our main results, where we comment on the empirical path measure interpretation of Bose-Einstein condensation at the end of that section. Section~\ref{sec-proofs} is devoted to the proofs, and in the appendices in Section~\ref{sec-app} we provide the notion of large deviations principle, and review facts on the zeta function of Riemann and on quantum statistical mechanics for systems of Bosons.

\section{Results}\label{sec-Results}
In this section we present our results. The phase transition result in Theorem~\ref{THManalysisRF} leads us to a probabilistic interpretation of Bose-Einstein condensation on the path measure level. I
Throughout the paper, we fix $ \beta >0 $.  Let $ \Omega:=\{\omega\colon [0,\infty)\to\R^d\colon \omega\mbox{ continuous}\} $ be the set of continuous functions $ [0,\infty)\to \R^d $. The elements in $ \Omega $ are called trajectories or paths and we denote by $ \Omega_k=\{\omega\colon [0,k\beta]\to\R^d\colon\omega\mbox{ continuous}\}, k\in\N $, the set of paths for time horizon $ [0,k\beta] $. We write $ \Omega_\beta $ for $ \Omega_1 $. We equip $\Omega $ (respectively $ \Omega_k $) with the topology of uniform convergence and with the corresponding Borel $\sigma$-field $ \Bcal $ (respectively $ \Bcal_k $). 
We consider $N$ Brownian motions, $ B^{\ssup{1}},\ldots, B^{\ssup{N}}$, with time horizon $ [0,\beta] $ as $ N $ random variables taking values in $\Omega_\beta$. For the reader's convenience, we repeat the definition of a Brownian bridge measure; see the Appendix in \cite{Sz98}. We decided to work with Brownian motions having generator $\Delta$ instead of $\frac 12\Delta$. We write $\P_x$ for the probability measure under which $B=B^{\ssup{1}}$ starts from $x\in\R^d$. The canonical (non-normalised) Brownian bridge measure on the time interval $[0,\beta]$ with initial site $x\in\R^d$ and terminal site $y\in \R^d$ is defined as
$ \mu_{x,y}^\beta(A)=\P_x(B\in A;B_\beta\in\d y)/\d y $ for $ A\subset\Omega_\beta $ measurable. Hence, the Brownian bridge measure for a Brownian bridge confined to a subset $ \L_N\subset \R^d $ is defined by
\begin{equation}\label{nnBBM}
\mu_{x,y}^{\beta,N}(A)=\frac {\P_x(B\in A;B_\beta\in\d y,B_{[0,\beta]}\subset\L_N)}{\d y},\qquad A\subset \Omega_\beta\mbox{ measurable.}
\end{equation}
Then $\mu_{x,y}^{\beta,N} $ is a regular Borel measure on $\Omega_\beta$.
The normalised Brownian bridge measure is defined as 
\begin{equation}\label{def-normalisedBBm}
\P_{x,y}^{\beta,N}=\mu_{x,y}^{\beta,N}/\mu_{x,y}^{\beta,N}(\Omega_\beta), 
\end{equation} which is a probability measure on $ \Omega_\beta $. 

We conceive the empirical path measure as a random element in $ \Mcal_1(\Omega) $, hence,  we need a convenient extension of any continuous path $ [0,\beta]\to\R^d $ to a continuous path $ [0,\infty)\to\R^d $ in the definition \eqref{delemp} of the empirical path measure. 
For any $ x\in\R^d $ we denote by $ P^x $ the Brownian probability measure on $ \Omega$, i.e., the canonical Wiener measure with deterministic start in $ x\in\R^d $ (\cite{Gin70}). In the following we write alternatively $ \omega_t $ or $ \omega(t) $ for any point of a path $ \omega $. 
Given a path $ \omega\in\Omega_\beta $ with time horizon $ [0,\beta] $ define
\begin{equation}\label{defprob}
P^{\ssup{\beta}}_{\omega}=\delta_{\omega}\otimes_\beta P^{\omega_\beta(\beta)}\in\Mcal_1(\Omega,\Bcal),
\end{equation}
where the product $ \otimes_\beta $ is defined for the ''splice`` of two paths, i.e., for $ \omega\in\Omega_\beta $ and $ \widetilde \omega\in\Omega $ define $ \overline{\omega}\in\Omega $ by $ \overline{\omega}(t)=\omega(t\wedge\beta), t\in [0,\infty) $, and $ \omega\otimes_\beta\widetilde \omega\in\Omega $ such that $ \omega\otimes_\beta\widetilde \omega=\widetilde \omega $ if $ \widetilde \omega(0)\not=\omega(\beta) $ and 
\begin{equation}\label{defsplice}
\omega\otimes_\beta\widetilde \omega(t)=\left\{\begin{array}{r@{\; \;}l}
\omega(t) & \mbox{ for }t\in [0,\beta]\\
\widetilde \omega (t-\beta) & \mbox{ for }t\in (\beta,\infty)
                                                     \end{array}\right.
\end{equation}
if $ \widetilde \omega(0)=\omega(\beta) $. The mapping $ \omega\in\Omega_\beta\mapsto P^{\ssup{\beta}}_{\omega}\in\Mcal_1(\Omega,\Bcal) $ is measurable, and the family $ \{P^{\ssup{\beta}}_{\omega}\colon \omega\in\Omega_\beta\} $ satisfies the Markov property (\cite[Lemma~4.4.21]{DS01}).
Hence, the empirical path measure 
\begin{equation}\label{delemp}
\begin{aligned}
L_N\colon\Omega_\beta^N\to \Mcal_1(\Omega), \omega\mapsto L_N(\omega)=\frac{1}{N}\sum_{i=1}^N\delta_{\omega^{\ssup{i}}}\otimes_\beta P^{\omega^{\ssup{i}}_\beta},
\end{aligned}
\end{equation} is $ \Omega_\beta^{\otimes N} $ measurable. Here $ \omega=(\omega^{\ssup{1}},\ldots,\omega^{\ssup{N}})\in\Omega_\beta^N $. Our main result concerns a large deviations principle for the distributions of $ L_N $ under the symmetrised measure $ \P_N^{\ssup{\rm sym}} $. Recall that $ \P_N^{\ssup{\rm sym}} $ is a probability measure on $ \Omega_\beta^{\otimes N} $.

To formulate the rate functions we need some notations.

Now we introduce the rate function.
Let 
$$ 
\Mcal=\{Q\in [0,1]^\N\colon \sum_{l\in\N}Q(l)=1, Q(l)\ge Q(l+1)\,\forall l\in\N\} 
$$ be the set of monotonously non-increasing probability functions on $ \N $. For $ Q\in\Mcal $ define $  \widehat Q(k)=Q(k)-Q(k+1) $ for any $ k\in\N $. For $ d\ge 1 $ let
\begin{equation}\label{criticalQ}
\widehat Q^*(k)=\frac{1}{\rho(4\pi\beta)^{d/2}k^{1+\d/2}},\quad k\in\N,
\end{equation} be given, and define the functional

\begin{equation}\label{functionalS}
\Scal(Q)=\sum_{k=1}^\infty\widehat Q(k)\Big(\log\frac{\widehat Q(k)}{\widehat Q^*(k)}-1\Big),\quad Q\in\Mcal.
\end{equation}

The rate function is given by

\begin{equation}\label{defratefctsym}
I^{\ssup{\rm sym}}(\mu)=\inf_{Q\in\Mcal} \Big\{S(Q)+I^{\ssup{Q}}(\mu)
\Big\}-\chi(\beta,\rho),\quad \mu\in\Mcal_1(\Omega),
\end{equation}

where

\begin{equation}\label{defrateQ}
I^{\ssup{Q}}(\mu)=\sup_{F\in\Ccal_{\rm b}(\Omega)}\Big\{\langle F,\mu\rangle-\sum_{k\in\N}\widehat Q(k)\log\E_{0,0}^{k\beta}\Big({\rm e}^{F(B)}\Big)\Big\},\quad \mu\in\Mcal_1(\Omega),
\end{equation}
and where the function $ \chi(\beta,\rho):=\inf_{Q\in\Mcal}\{\Scal(Q)\}$ is given as the negative logarithmic limit of the partition function $ Z_N^{\ssup{\rm sym}}(\beta) $, see Theorem~\ref{THMfree}, and where $ \Ccal_{\rm b}(\Omega) $ is the space of continuous bounded functions of the paths in $ \Omega $. $ \E_{0,0}^{k\beta} $ denotes the expectation with respect to the Brownian bridge probability measure $ \P_{0,0}^{k\beta} $ defined in \eqref{def-normalisedBBm} extended as in \eqref{defprob} to a probability measure in $ \Mcal_1(\Omega) $. Here, $ I^{\ssup{Q}} $ is a Fenchel-Legendre transform, but {\it not\/} the one of a logarithmic moment generating function of any random variable.
In particular, $ I^{\ssup{Q}} $, and therefore also $ I^{\ssup{\rm sym}} $, are nonnegative, and $ I^{\ssup{Q}} $ is convex as a supremum of linear functions. There seems to be no way to represent $ I^{\ssup{Q}}(\mu) $ as the relative entropy of $ \mu $ with respect to any measure.

We need the following limit for the partition function to derive the large deviations principle.

\begin{theorem}[{\bf Limit of the partition function}]\label{THMfree} Let $ \rho\in(0,\infty) $ and $ \L_N\subset\R^d $ with $ \L_N\uparrow\R^d $ and $ N/|\L_N|\to\rho $ as $ N\to\infty $ and with $ \dist(0,\partial\L_N)=cN $ for some $ c>0 $. Then
\begin{equation}\label{variationalfreeenergy}
\lim_{N\to\infty}\frac{1}{N}\log Z_N^{\ssup{\rm sym}}(\beta)=-\inf_{Q\in\Mcal}\Big\{\Scal(Q)\Big\}.
\end{equation}
\end{theorem}

This Theorem gives the normalisation part $ \chi(\beta,\rho)=\inf_{Q\in\Mcal}\{\Scal(Q)\} $ in the definition of the rate function $ I^{\ssup{\rm sym}} $ in \eqref{defratefctsym} for the large deviations principle in the following Theorem~\ref{THMLDP-LN}. 

\begin{theorem}[{\bf Large deviations for $ \mathbf{L_N} $}]\label{THMLDP-LN}
Let $ \rho\in(0,\infty) $ and $ \L_N\subset\R^d $ with $ \L_N\uparrow\R^d $ and $ N/|\L_N|\to\rho $ as $ N\to\infty $ and with $ \dist(0,\partial\L_N)=cN $ for some $ c>0 $.

Under the symmetrised measure $\P_N^{\ssup{\rm sym}}$ the empirical path measures $ (L_N)_{N\in\N} $ satisfy a large deviations principle on $ \Mcal_1(\Omega) $ with speed $ N $ and rate function $ I^{\ssup{\rm sym}} $.
\end{theorem}

The proof of Theorem~\ref{THMLDP-LN} is in Subsection~\ref{proof-secLDP-LN}. The proof is based on the following variational formula for the logarithmic limit of the normalisation in Theorem~\ref{THMfree} and on the analysis of the rate function in Theorem~\ref{THM-variational} below. 

We give a brief informal interpretation of the shape of the rate functions in \eqref{defratefctsym} and \eqref{defrateQ}. As remarked earlier, the symmetrised measure $ \P_N^{\ssup{{\rm sym}}} $ arises from a two-step probability mechanism. This is reflected in the representation of the rate function $ I^{\ssup{{\rm sym}}}$ in \eqref{defratefctsym}: in a peculiar way (we describe it at the beginning of Section~\ref{sec-proofs}), the term $ S(Q)-\chi(\beta,\rho) $ describes the large deviations of the discrete empirical shape measure for integer partitions. The discrete empirical shape measure $ Q_N $, defined as
\begin{equation}\label{shapemeasure}
Q_N\colon\Pcal_N\to\Mcal_1(\N), \l\mapsto\frac{1}{N}\sum_{k=\cdot}^Nr_k(\l), 
\end{equation} governs a particular distribution of $N$ independent, but not identically distributed, Brownian bridges. Under this distribution, $L_N$ satisfies a large deviations principle with rate function $I^{\ssup{Q}}$, which can also be guessed from the G\"artner-Ellis theorem \cite[Th.~4.5.20]{DZ98}. The presence of a two-step mechanism makes it impossible to apply this theorem directly to $ \P_{N}^{\ssup{{\rm sym}}} $.

Let us contrast this to the case of i.i.d.~Brownian bridges $B^{\ssup 1},\dots,B^{\ssup N}$, starting in the origin, i.e., we replace $ \P_{N}^{\ssup{{\rm sym}}} $  by $ (\P_{0,0}^\beta)^{\otimes N}$. Here the empirical path measure $ L_N $ satisfies a large deviations principle with rate function 
$$
I(\mu)=\sup_{F \in\Ccal_{\rm b}(\Omega)}\Bigl\{\langle F ,\mu\rangle-\log\E_{0,0}^\beta\bigl({\rm e}^{F(B)}\bigr)\Bigr\},
$$
as follows from an application of Cram\'er's theorem \cite[Theorem~6.1.3]{DZ98}.  Note that $I(\mu)$ is the relative entropy of $\mu$ with respect to $ \P_{0,0}^\beta\circ B^{-1}$. Although there is apparently no reason to expect a direct comparison between the distributions of $L_N$ under  $ \P_{N}^{\ssup{{\rm sym}}} $ and under $ (\P_{0,0}^\beta)^{\otimes N}$, the rate functions admit a simple relation: it is easy to see that $ I^{\ssup{Q}}\ge I $ for the measure $ Q\in\Mcal $ with $ \widehat Q(k)=\delta_1 $, since
\begin{eqnarray}
\begin{aligned}
-\sum_ {k=1}^\infty\widehat Q(k)\log\E^{k\beta}_{0,0}\big[{\rm e}^{F(B)}\big]&\ge -\log\E_{0,0}^\beta\bigl({\rm e}^{F(B)}\bigr). 
\end{aligned}
\end{eqnarray} 
In particular, $ I^{\ssup{\rm sym}}\ge I $.

\begin{remark}
The techniques of the proof of Theorem~\ref{THMLDP-LN} apply also to a proof of a large deviations principle under the symmetrised measure $ \P_N^{\ssup{\rm sym}} $ for the empirical path measure $ \widetilde L_N=1/N\sum_{i=1}^N\delta_{B^{\ssup{i}}} $, which is a random element in $ \Mcal_1(\Omega_\beta) $. The rate function is
$$
\widetilde I^{\ssup{\rm sym}}(\mu)=\inf_{Q\in\Mcal}\Big\{\Scal(Q)-\sup_{F\in\Ccal_{\rm b}(\Omega_\beta)}\Big\{\langle F,\mu\rangle-\sum_{k=1}^\infty\widehat Q(k)\log\E_{0,0}^{k\beta}\Big({\rm e}^{\sum_{l=0}^{k-1}F(B_{[l\beta,(l+1)\beta]})}\Big)\Big\}\Big\}\quad, \mu\in\Mcal_1(\Omega_\beta).
$$
Similar results hold for the mean of the occupation $ Y_N $ measures,
$$
Y_N=\frac{1}{N}\sum_{i=1}^N\frac{1}{\beta}\int_0^\beta\delta_{B^{\ssup{i}}_s}\d s\in\Mcal_1(\R^d).
$$
However, these rate functions seem not to give enough information to derive the phase transition as in Theorem~\ref{THManalysisRF}, and to obtain a probabilistic interpretation of Bose-Einstein condensation.
\end{remark}
Our large deviations result is accompanied by an analysis of the variational formula for the rate function $ I^{\ssup{\rm sym}} $ \eqref{defratefctsym}, i.e., the analysis for zeros of the rate function. This gives the proof phase transition for empirical path measures depending on the dimension and the density parameter in Theorem~\ref{THManalysisRF}. We analyse first the variational formula for the limit of the partition function.
\begin{remark}[{\bf Free energy}]
The variational formula for the limit in \eqref{variationalfreeenergy} gives the {\bf specific free energy}
$ f(\beta,\rho):=\lim_{N\to\infty}-1/\beta|\L_N|Z_N^{\ssup{\rm sym}}(\beta) $ for  inverse temperature $ \beta $ and density $ \rho $, i.e.,
\begin{equation}
f(\beta,\rho)=\frac{\rho}{\beta}\inf_{Q\in\Mcal}\Big\{\sum_{k=1}^\infty\widehat Q(k)\log \Big(\frac{\widehat Q(k)}{\widehat Q^*(k)}-1\Big)\Big\}.
\end{equation}
This is the specific free energy for an infinite systems of non-interacting Bosons in the thermodynamic limit, i.e., the limit where $ N/|\L_N|\to\rho $ as $ N\to\infty $. 
\end{remark}

We analyse the variational formula for $ \chi(\beta,\rho) $ \eqref{variationalfreeenergy}, and we derive an expression for the free energy $ f $ as a function of $ \beta $ and $ \rho $. Define a dimension dependent critical density
\begin{equation}
\rho_c=\left\{\begin{array}{r@{\;,\;}l}\frac{1}{(4\pi\beta)^{d/2}}\zeta\Big(\frac{d}{2}\Big) & \mbox{ for } d\ge 3\\
+\infty & \mbox{ for } d=1,2
              \end{array}\right.,
\end{equation}
where $ \zeta $ is the zeta function of Riemann,  
$$
\zeta\Big(\frac{d}{2}\Big)=\sum_{k=1}^\infty k^{-\frac{d}{2}}.
$$
Furthermore, denote by $ g_{s}(\alpha) $ the so-called Bose functions (see \eqref{defBosefunctionapa} in Appendix~\ref{appendix-Bosefunction})
$$
g_{s}(\alpha)=\sum_{k=1}^\infty k^{-s}{\rm e}^{-\alpha k}\quad\mbox{ for all } \alpha>0\mbox{ and all } s>0. 
$$
For any $ \rho<\rho_c $ we denote by $ \alpha=\alpha(\beta,\rho) $ the unique root of
\begin{equation}\label{root}
\rho=\frac{1}{(4\pi\beta)^{d/2}}\sum_{k=1}^\infty k^{-d/2}{\rm e}^{-\alpha k}.
\end{equation}
The essential difference in $ d\ge 3 $ and $ d=1,2 $ lies in the fact that in the latter two cases the corresponding Bose functions, $ g_1(\alpha) $ respectively $ g_{\frac{1}{2}}(\alpha) $, diverge as $ \alpha\to 0 $ (see Appendix~\ref{appendix-Bosefunction} and \cite{Gram25}). For  $ d=1,2 $ there is a unique $ \alpha $ for any density $ \rho <\infty $. For $ d\ge 3 $ there is such an unique $ \alpha $ given only for densities $ \rho <\rho_c $. Hence, this is the mathematical origin of the so-called Bose-condensation (see Appendix~\ref{qm-sec}), where for $ d\ge 3 $ and $ \rho>\rho_c $ particles condense in the zero mode state. This corresponds to the phase transition for the empirical path measure, when the empirical path measure has positive mass on paths which do not return to the origin at any times which are an integer  multiple of $ \beta $.

In the following theorem we analyse the variational problem \eqref{variationalfreeenergy}.

\begin{theorem}[{\bf Analysis of the variational formula for $ \chi(\beta,\rho)$}]\label{THM-variational}
For any $ \rho<\infty $ in dimensions $ d=1,2 $, and  $ \rho<\rho_c $ in dimensions $ d\ge 3 $, there is a unique minimiser $ Q\in\Mcal $ of the variational formula \eqref{variationalfreeenergy} with
\begin{equation}
\widehat Q(k)=\frac{{\rm e}^{-\alpha k}}{\rho(4\pi\beta)^{d/2}k^{1+\frac{d}{2}}}\quad\mbox{ for } k\in\N,
\end{equation}
whereas for dimensions $ d\ge 3 $ and densities $ \rho>\rho_c $, there is no minimiser for the variational problem \eqref{variationalfreeenergy}, but the infimum is attained for any minimising sequence $ (Q_n)_{n\in\N} $ of $ Q_n\in\Mcal $ such that $ Q_n\to Q^* $ as $ n\to\infty $.

The specific free energy for $ d\ge 3 $ is given by
\begin{equation}\label{freeed3}
f(\beta,\rho)=\left\{\begin{array}{r@{\;,\;}l} -\frac{1}{(4\pi\beta)^{d/2}\beta}g_{\frac{d+2}{2}}(\alpha)-\frac{1}{\beta}\rho\alpha & \mbox{ for } \rho<\rho_c\\[1.5ex]
-\frac{1}{(4\pi\beta)^{d/2}\beta}\zeta\Big(\frac{d+2}{2}\Big) & \mbox{ for } \rho>\rho_c,
                     \end{array}\right.
\end{equation}
and for $ d=1,2 $ by
\begin{equation}\label{freeed12}
f(\beta,\rho)=-\frac{1}{(4\pi\beta)^{d/2}\beta}g_{\frac{d+2}{2}}(\alpha)-\frac{\rho\alpha}{\beta},
\end{equation}
where $ \alpha $ is the unique root of \eqref{root}.
\end{theorem}

This result leads to the analysis of the rate function $ I^{\ssup{\rm sym}} $. Let
\begin{equation}\label{defAK}
A_k=\{\omega\otimes_{k\beta}\xi\colon\omega\in\Omega_k,\omega(0)=\omega(k\beta),\xi\in\Omega\}\subset\Omega, k\in\N,
\end{equation}
be the set of paths in $ \Omega $ which result from the splice \eqref{defsplice} of Brownian bridges paths of time horizon $ [0,k\beta] $ with any path $ \xi\in\Omega $.

\begin{theorem}[{\bf Analysis of the rate function $ I^{\ssup{\rm sym}}$}]\label{THManalysisRF}
Under the assumptions of Theorem~\ref{THMLDP-LN} 
the following holds:
\begin{enumerate}
\item[(i)] $ d=1,2 $.  A unique minimiser $ \mu^*\in\Mcal_1(\Omega) $ of the rate function $ \mu\mapsto I^{\ssup{\rm sym}}(\mu) $ is given with $ \sum_{k\in\N}k\mu^*(A_k)=1 $.

\item[(ii)] $ d\ge 3 $ and $ \rho<\rho_c $.  A unique minimiser $ \mu^*\in\Mcal_1(\Omega) $ of the rate function $ \mu\mapsto I^{\ssup{\rm sym}}(\mu) $ is given with $ \sum_{k\in\N}k\mu^*(A_k)=1 $. 

\noindent For $ \rho>\rho_c $ there is no unique minimiser be given, but there exist minimising sequences $ (\mu_n)_{n\ge 1}, \mu_n\in\Mcal_1(\Omega) $,  with $ \sum_{n=1}^\infty k\mu_n(A_k)=1 $ for any $ n\in\N $ such that $ \mu_n\to\mu^0\in\Mcal_1(\Omega) $ weakly as $ n\to\infty $ with $ \sum_{n=1}^\infty k\mu^0(A_k)<1$.
\end{enumerate}
\end{theorem}

\begin{proofsect}{Proof}
$ \Scal_\chi:=\Scal(Q)-\chi(\beta,\rho) \ge 0 $ for all $ Q\in\Mcal $. Moreover, $ I^{\ssup{Q}}(\mu)\ge 0 $ for any $ Q\in\Mcal $, which can be seen putting $ F\equiv 0 $ in the supremum \eqref{defrateQ}. The functional $ \Scal $ can be written as
$$
\Scal(Q)=qH(\widehat P|\widehat P^*)+q\log\frac{q}{q^*}-q\quad, Q\in\Mcal, 
$$ where $ q=\sum_{k=1}^\infty\widehat Q(k) , q*=\sum_{k=1}^\infty\widehat Q^*(k) $, and $ H $ is the relative entropy of the two probability measures $ \widehat P=q^{-1}\widehat Q $ and $ \widehat P^*=(q^*)^{-1} \widehat Q^* $. 
Hence, the level sets of $ Q\mapsto \Scal_\chi(Q) $ are compact.

\noindent (i) Case $ d=1,2 $.
Assume that $ \mu $ is a zero of $ I^{\ssup{\rm sym}} $. Since the level sets are compact, there is a $ Q^0\in\Mcal $ that minimises the formula on the right hand side in the definition of $ I^{\ssup{\rm sym}} $ \eqref{defratefctsym}. Clearly $ I^{\ssup{Q^0}}<\infty $. As both parts of the rate function $ I^{\ssup{\rm sym}} $ are positive, we conclude that $ Q^0 $ is the unique minimiser of the functional $ \Scal $ (see Subsection~\ref{proof-sec-variational} in the proof of Theorem~\ref{THM-variational}), i.e., $ \Scal_\chi(Q^0)=0 $. Note further that a subclass of $ \Ccal_{\rm b}(\Omega) $ consists of those function $ F_k, k\in\N $,  such that 
$$ 
F_k(\omega\otimes_{k\beta}\xi)=\int_0^{k\beta}f(\omega_s)\d s\quad \mbox{ for all }\xi\in\Omega,
$$
for continuous bounded functions $ f_k\colon\R^d\to R $
The Euler-Lagrange equations yield the following system of equations
$$
\langle H,\mu\rangle=\sum_{k\in\N}\widehat Q^0(k)\frac{\E_{0,0}^{k\beta}\Big(H{\rm e}^{F(B\otimes_{k\beta}P^0)}\Big)}{\E_{0,0}^{k\beta}\Big({\rm e}^{F(B\otimes_{k\beta}P^0)}\Big)},
$$
where $ H = F_k , k\in\N $, is any of the functions defined above. $ F\equiv 0 $ solves this system, and if one chooses $ f_k=1/\beta ,k\in\N $, in the definition of the functional $ F_k $, one obtains that
$$
k\mu(A_k)=k\widehat Q^0(k)\quad, k\in\N. 
$$ 
This identifies $ \mu $ as $ \mu^* $, since $ \sum_{k=1}^\infty k\widehat Q^0(k)=\sum_{k=1}^\infty Q(k)=1 $. 

\noindent (ii) Case $ d\ge 3 $. For densities $ \rho<\rho_c $ the assertion follows as in (i). 
Let $ \rho>\rho_c $.  Define the sequence $ (Q_n)_{n\in\N} , Q_n\in\Mcal $, by
\begin{equation}\label{defseq}
\widehat Q_n(k)=\left\{\begin{array}{r@{\;,\;}l}
\widehat Q^*(k) & \mbox{ for } k\not= n\\
\widehat Q^*(k)+\frac{(\rho-\rho_c)}{n\rho} & \mbox{ for } k=n
                       \end{array}\right.,
\end{equation} 
where $ \widehat Q^* $ is given in \eqref{criticalQ} (see Subsection~\ref{proof-sec-free}), i.e., 
$$
\widehat Q^*(k)=\frac{1}{\rho(4\pi\beta)^{d/2}k^{1+\d/2}},\quad k\in\N.
$$
Clearly, $ Q_n\to Q^* $ strongly as $ n\to\infty $, and $\lim_{n\to\infty}\Scal_\chi(Q_n)=0 $. 
As 
$$
\sum_{k=1}^\infty Q^*(k)=\sum_{k=1}^\infty k\widehat Q^*(k)=\frac{1}{\rho(4\pi\beta)^{d/2}}\zeta(d/2)=\frac{1}{\rho(4\pi\beta)^{d/2}}g_{d/2}(0) <1,
$$
the infimum is not attained in $ \Mcal $. Hence, from the considerations in (i), we get a sequence $ (\mu_n)_{n\in\N} $ of zeros $ \mu_n\in\Mcal_1(\Omega) $ of the rate function $ I^{\ssup{Q_n}} $ with $ \mu_n(A_k)=\widehat Q_n(k) $ and 
$ \sum_{k=1}^\infty k\mu_n(A_k)=1 $ such that $ I^{\ssup{Q_n}}(\mu_n)=0 $. 
Then,
$$
I^{\ssup{\rm sym}}(\mu_n)\le\inf_{Q\in\Mcal}\Big\{\Scal_\chi(Q)+I^{\ssup{Q}}(\mu_n)\Big\}\le \Scal_\chi(Q_n),
$$
and from the upper semi continuity of the functional $ \Scal $, the strong convergence of $ Q_n\to Q^* $ as $ n\to\infty $, and from the fact that the level sets of $ I^{\ssup{\rm sym}} $ are compact, one derives that $ \mu_n\to\mu^0\in\Mcal_1(\Omega) $ weakly as $ n\to\infty $, and $ \lim_{n\to\infty} I^{\ssup{\rm sym}}(\mu_n)=0 $.  
\qed
\end{proofsect}

Let us draw an easy corollary from this theorem.

\begin{cor}[{\bf Law of large numbers}]
Under the assumptions of Theorem~\ref{THManalysisRF} the following holds.
\begin{enumerate}
\item[(i)] For $ d=1,2 $, and any density $ \rho<\infty $, there is a law of large numbers. Under the probability measure $ \P_N^{\ssup{\rm sym}} $, the sequence $ (L_N)_{N\in\N} $ converges in distribution to the measure $ \mu^*\in\Mcal_1(\Omega) $. 
\item[(ii)] For $ d\ge 3 $ and $ \rho<\rho_c $ there is a law of large numbers. Under the probability measure $ \P_N^{\ssup{\rm sym}} $, the sequence $ (L_N)_{N\in\N} $ converges in distribution to the measure $ \mu^*\in\Mcal_1(\Omega) $. 
\end{enumerate}
\end{cor}

The main conclusion of the large deviations principle in Theorem~\ref{THMLDP-LN} and Theorem~\ref{THManalysisRF} is the following phase transition for the mean empirical path measure, which gives a path measure interpretation of Bose-Einstein condensation (BEC).

\noindent {\bf Path measures and their interpretation as Bose-Einstein condensation}

Let $ N $ Brownian motions with time horizon $ [0,\beta] $ confined in sets $ \L_N\subset \R^d $ given such that $ \L_N\uparrow\R^d $ and $ \dist(0,\partial\L_N)=cN $ for some $ c>0 $ and $ N/|\L_N|\to\rho\in (0,\infty) $ as $ N\to\infty $. Then the following holds: 
\begin{enumerate}

\item[(i)] For $ \beta> 0 $ there is a $ \rho_c=\rho_c(\beta,d) $ such that:

\noindent {\bf no BEC:} Case $ \rho<\rho_c$ for $ d\ge 3 $, $ \rho>0 $ for $ d=1,2 $:

\noindent  $ L_N\to\mu^*\in\Mcal_1(\Omega) $ under $ \P_N^{\ssup{\rm sym}} $ as $ N\to\infty $ with $ \sum_{k=1}^\infty k\mu^*(A_k)=1 $ 

\noindent {\bf BEC:} Case $ \rho<\rho_c $ and $ d\ge 3 $: 

\noindent $ L_N\to\mu^0\in\Mcal_1(\Omega) $ under $ \P_N^{\ssup{\rm sym}} $ as $ N\to\infty $ with $ \sum_{k=1}^\infty k\mu^0(A_k)<1 $.

\item[(ii)] For $ \rho\in(0,\infty) $ there exists a 
$$ 
\beta_c=\left\{\begin{array}{r@{\;,\;}l}
\frac {1}{4\pi}\big(\frac{\rho}{\zeta(d/2)}\Big)^{2/d} & \mbox{ for } d\ge 3\\
+\infty & \mbox{ for } d=1,2
               \end{array}\right.,
$$ 
such that: 

\noindent {\bf no BEC:} Case $ \beta<\beta_c$ for $ d\ge 3 $ and $ \beta>0 $ for $ d=1,2 $:

\noindent  $ L_N\to\mu^*\in\Mcal_1(\Omega) $ under $ \P_N^{\ssup{\rm sym}} $ as $ N\to\infty $ with $ \sum_{k=1}^\infty k\mu^*(A_k)=1 $ 

\noindent {\bf BEC:} Case $ \beta>\beta_c$ and $ d\ge 3$:

\noindent  $ L_N\to\mu^0\in\Mcal_1(\Omega) $ under $ \P_N^{\ssup{\rm sym}} $ as $ N\to\infty $ with $ \sum_{k=1}^\infty k\mu^0(A_k)<1 $.
\end{enumerate}

If $ d=1,2 $, or $ \rho<\rho_c $ for $ d\ge 3 $, the mean empirical path measure has support on those paths in which one can insert, starting from time origin, a concatenation of any finite number of Brownian motions with time horizon $ [0,\beta] $, i.e., for any $ k\in\N $ one can find in paths $ \omega_k\in A_k $ exactly $ k $ Brownian motions concatenated to a Brownian bridge with horizon $ [0,k\beta] $. This follows from the concatenation of the Brownian motions due to the cycle structure of the permutations and due to the Lebesgue integration of any initial position in the definition of the symmetrised measure $ \P_N^{\ssup{\rm sym}} $ \eqref{symPdef1}. If the density $ \rho $ is high enough for $ d\ge 3 $, i.e., $ \rho>\rho_c $ (or equivalently, if the inverse temperature is sufficiently large for given density, i.e., $ \beta>\beta_c $, for $ d\ge 3 $), the mean path measure has positive weight for paths with an infinite time horizon, that is, concatenation of  any finite number of Brownian motions with time horizon $ [0,\beta] $, i.e., any finite cycle path in $ A_k $, is not sufficient, because there is an excess density $ (\rho-\rho_c) $ of Brownian motions with time horizon $ [0,\beta] $. These motions  concatenate to infinite long cycle, that is, these cycles grow with the system size in the thermodynamic limit.  The fraction of these motions is 
\begin{equation}
1-\frac{\rho_c}{\rho}=1-\Big(\frac{\beta_c}{\beta}\Big)^{d/2}.
\end{equation}

Let us briefly contrast our results to the ones derived in \cite{AD06} and \cite{AK06}. There the empirical path measure was studied under the symmetrised measure
$$
\P_{m,N}^{\ssup{{\rm sym}}}=\frac{1}{N!}\sum_{\s\in\Sym_N}\;\int_{\R^d}\cdots \int_{\R^d} m(\d x_1)\cdots m(\d x_N)\bigotimes_{i=1}^N\P^{\beta}_{x_i,x_{\s(i)}},
$$
where $ m\in\Mcal_1(\R^d) $ is an initial distribution for the Brownian motions. Multilevel large deviations principle are obtained, however, a phase transition for the mean path measure does not exist. Only in the special case where $ m $ is the Lebesgue measure restricted to a fixed box $ \L\subset\R^d $, it turns out that the rate function for finite time horizon in the limit $ N\to\infty $ is precisely the Donsker-Varadhan rate function. This rate function appears in the limit of infinite time horizon, hence its appearance is an indication of the presence of infinite long cycles. Note, this corresponds to the limit $ N/|\L|\to\infty $, i.e., infinite density. In \cite{AD06} a similar result was shown for quantum spin system with mean-field interaction.

Our results  here are essentially different, because of the definition of the measure $ \P_N^{\ssup{\rm sym}} $ and of the limit $ N/|\L_N|\to\rho\in(0,\infty) $, and go far beyond the ones in \cite{AD06} and \cite{AK06}. 

Finally, let us remark that our results show that the discrete empirical shape measure
$$
Q_N\colon\Pcal_N\to\Mcal_1(\N),\l\mapsto Q_N^\l(\cdot)=\frac{1}{N}\sum_{k=\cdot}^Nr_k(\l),
$$
shall satisfies a large deviations principle under the distributions $ \mu_N $ given by
$$
\mu_N(\l)=\frac{1}{Z_N^{\ssup{\rm sym}}}(\beta)\prod_{k=1}^N\Big(\frac{|\L_N|^ {r_k}}{r_k!k^{r_k}}\Big)\bigotimes\big(\mu_{0,0}^{k\beta}\big)^{r_k}(\Omega_\beta^N).
$$
The distribution $ \mu_N $ is not the uniform distribution on the partitions $ \Pcal_N $ ($d\ge 3$). Limit theorems for scaled shape measures have been obtained in \cite{Ver96} for uniform distributions of partitions. Large deviations under the uniform distribution from this limit behaviour  with appropriate scaling of the shape measures have been derived in \cite{DVZ00}. 
In \cite{A07b} we will study such large deviations for the non-uniform distributions given above.


\section{Proofs}\label{sec-proofs}
In this section we proof our the Theorems~\ref{THMfree},~\ref{THM-variational} and \ref{THMLDP-LN}. In Subsection~\ref{proof-sec-free} we prove the limit of the partition function, i.e., the thermodynamic limit of the free energy and its expression through a variational formula for discrete shape measures. 
The large deviations principle is proved in Subsection~\ref{proof-secLDP-LN}, and the analysis of the variational formula is in Subsection~\ref{proof-sec-variational}.

\subsection{Proof of Theorem~\ref{THMfree}}\label{proof-sec-free}
The key idea is to use the cycle representation for permutations and to replace the sum over permutations in the definition \eqref{partitionfunction} of the partition function by a sum over integer partitions of $ N $. 
For any partition $ \lambda\in\Pcal_N $ we assign the discrete empirical shape measure $ Q_N^\lambda $ for the  occupation numbers  $ \{r_k\}_{k=1}^N $ with $ \sum_{k=1}^Nk r_k=N $, i.e., the probability measure $ Q^\l_N\in\Mcal_1(\N) $, defined as
\begin{equation}\label{discreteshapemeasure}
\begin{aligned}
Q_N&\colon\Pcal_N\to\Mcal_1(\N)\\
&\lambda\mapsto Q_N^\lambda(\cdot)=\frac{1}{N}\sum_{k=\cdot}^Nr_k(\lambda).
\end{aligned}
\end{equation}
$ Q_N^\lambda(l)=0 $ for $ l>N $. 
Here $ Q_N^\lambda(1)=\sum_{k=1}^Nr_k $ is the number of components/cycles. The mapping $ Q_N $ is injective, which can easily be seen as follows. $ Q_N^\lambda=Q_N^{\lambda^\prime} $ for $\lambda,\lambda^\prime\in\Pcal_N $ with $ \lambda\not=\lambda^\prime $ implies successively $ r_N(\lambda)=r_N(\lambda^\prime),\ldots r_1(\lambda)=r_1(\lambda^\prime) $ and hence $ \lambda=\lambda^\prime $ in contradiction with our assumption.
For any partition $ \lambda\in\Pcal_N $, i.e., for any set of occupation numbers $ \{r_k\}_{k=1}^N $ with $ \sum_{k=1}^Nk r_k=N $, we have
$ N!/\prod_{k=1}^N r_k! k^{r_k} $ permutations of type $ \{r_k\}_{k=1}^N $. For a given permutation $ \l\in\Pcal_N $ we regroup the product of the Brownian bridges as follows. In the following we write $ r_k $ for $ r_k(\l) $. For any $ k\in\{1,\ldots,N\} $ with $ r_k\ge 1 $ we have exactly $ r_k $ groups/cycles of length $ k $ which involve altogether $ kr_k $ Brownian motions,
$$
\begin{aligned}
\Big(\int_{\L_N^k}\d x_1\cdots\d x_k\,\mu^{\beta,N}_{x_1,x_2}&\otimes\mu^{\beta,N}_{x_2,x_3}\otimes\cdots\mu^{\beta,N}_{x_k,x_1}\Big)\Big(\int_{\L_N^k}\d x_{k+1}\cdots\d x_{2k}\,\mu^{\beta,N}_{x_{k+1},x_{k+2}}\otimes\cdots\otimes\mu^{\beta,N}_{x_k,x_1}\Big)\\
&\cdots \Big(\int_{\L_N^k}\d x_1\cdots\d x_k\,\mu^{\beta,N}_{x_{(r_k-1)k+1},x_{(r_k-1)k+2}}\otimes\cdots\mu^{\beta,N}_{x_{r_k k},x_{(r_k-1)k+1}}\Big),
\end{aligned}
$$
numbered here from $ 1 $ to $ kr_k $. The structure of the product measure in each of the $ r_k $ groups/cycles is identical, hence we define the measure
\begin{equation}\label{defmuk}
\begin{aligned}
\mu_k&=\frac{1}{|\L_N|}\int_{\L_N}\d x\int_{\L_N^{k-1}}\d x_1\cdots\d x_{k-1}\,\mu_{x,x_1}^{\beta,N}\otimes\mu_{x_1,x_2}^{\beta,N}\otimes\cdots\otimes\mu_{x_{k-1},x}^{\beta,N}\\&=\int_{\L_N^{k-1}}\d x_1\cdots\d x_{k-1}\,\mu_{0,x_1}^{\beta,N}\otimes\mu_{x_1,x_2}^{\beta,N}\otimes\cdots\otimes\mu_{x_{k-1},0}^{\beta,N}=\mu_{0,0}^{k\beta,N},
\end{aligned}
\end{equation}
where we used spatial shift invariance of the Brownian motion and concatenated the Brownian bridges to a Brownian bridge of time length $ k\beta $ starting and terminating at the origin. Here, (see \eqref{nnBBM}),
$$
\mu_{0,0}^{k\beta,N}(A)=\frac{\P_0(A,B_{k\beta}\in\d y, B_{[0,k\beta]}\subset\L_N))}{\d y}\quad, A\subset\Omega_k \mbox{ measurable}.
$$
This concatenation is possible because of the indicator to stay within $ \L_N $. 
Hence, we obtain for the partition function
$$
\begin{aligned}
Z_N^{\ssup{\rm sym}}(\beta)&=\frac{1}{N!}\sum_{\s\in\Sym_N}\int_{\L_N}\d x_1\cdots\int_{\L_N}\d x_N\,\bigotimes_{i=1}^N\mu_{x_i,x_{\s(i)}}^{\beta,N}(\Omega_\beta^N)\\
&=\sum_{\l\in\Pcal_N}\prod_{k=1}^N\Big(\frac{|\L_N|^{r_k}}{(r_k)!k^{r_k}}\Big)\Big(\bigotimes_{k=1}^N(\mu_{0,0}^{k\beta,N})^{\otimes r_k}\Big)(\Omega_\beta^N),
\end{aligned}
$$


where we replaced the sum over permutations by a sum over partitions. To study the large $ N$-behaviour of the partition function we rewrite the expression on the right hand side of the latter expression. We sum over all elements in the image 
$$ 
\Mcal_N:=Q_N(\Pcal_N)= 
\subset\{Q\in \{0,1/N,\ldots,1\}^\N\colon \sum_{l=1}^NQ(l)= 1, Q(l)\ge Q(l+1)\,\forall l\in\N\}
$$ of the discrete empirical shape measure $ Q_N $, which are probability measures $ Q\colon\N\to\{0,1/N,\ldots,1\} $ with $ Q(l)\ge Q(l+1) $ for all $ l\in\N $ and $ \sum_{l=1}^NQ(l)= 1 $ and $ Q(l)=0 $ for all $ l>N $. The occupation numbers for a partition $ \l\in\Pcal_N $ are then given by $ r_k=N(Q_N(k)-Q_N(k+1)), k=1,\ldots,N $, where $ Q_N=Q_N^\l $. We write $ \widehat Q_N(k)=Q_N(k)-Q_N(k+1)$ for $ k=1,\ldots, N, $ in the following. 
Thus we obtain 
$$
\begin{aligned}
Z_N^{\ssup{\rm sym}}(\beta)=\sum_{Q\in\Mcal_N}\Big(\prod_{k=1}^N\frac{|\L_N|^{N\widehat Q(k)}}{(N\widehat Q(k))!k^{N\widehat Q(k)}}\Big)\Big(\bigotimes_{k=1}^N(\mu_{0,0}^{k\beta,N})^{\otimes N\widehat Q(k)}\Big)(\Omega_\beta^N). 
\end{aligned}
$$
On the right hand side we have to evaluate the probability mass on $ \Omega_\beta^N $. 
Now we use that $ \mu_{0,0}^{k\beta,N} , k\in\{1,\ldots,N\}, $ is a regular Borel measure on $ \Omega_k $ with total mass equal to 
\begin{equation}
\mu_{0,0}^{k\beta,N}(\Omega_k)=p_{k\beta,\L_N}(0,0).
\end{equation}
The transition density $ p_{k\beta,\L_N} $ is defined as
\begin{equation}\label{densitybox}
p_{k\beta,\L_N}(x,y)=\frac{\P_x(B_{k\beta}\in\d y)}{\d x}-r_{\L_N}(k\beta,x,y),\quad x,y\in \L_N,
\end{equation}
with
\begin{equation}
r_{\L_N}(k\beta,x,y)=\E_x(\tau_{\L_N}<k\beta;p_{k\beta-\tau_{\L_N}}(B_{\tau_{\L_N}},y)) 
\end{equation}
and $$ p_{k\beta}(x,y)=\frac{\P_x(B_{k\beta}\in\d y)}{\d x}=(4\pi\beta k)^{-d/2}{\rm e}^{-|x-y|^2/4\beta k}, $$ where $ \tau_{\L_N} $ is the exit time of the Brownian motion of the set $ \L_N $ (\cite{CZ95}). Without loss of generality we let $ \L_N \subset \R^d $ be centered such that $ \dist(0,\partial\L_N) =\alpha N $ for some $ \alpha>0 $. For any $ c $ with $ 0<k\beta<c^2/d $ and $ |x-y|<c\le \dist(x,\partial\L_N)\wedge\dist(y,\partial\L_N) $ one gets (\cite{CZ95}) the estimation
$$
r_{\L_N}(k\beta,x,y)\le (4\pi\beta k)^{-d/2}{\rm e}^{-\frac{c^2}{2k\beta}}. 
$$
Applying this with $ c=\sqrt{d\beta N} $ we get
$$
(4\pi\beta k)^{-d/2}(1-{\rm e}^{-dN/4\beta})\le \mu_{0,0}^{k\beta,N}(\Omega_k)\le (4\pi\beta k)^{-d/2}.
$$
Hence,
$$
\begin{aligned}
(1-{\rm e}^{-dN/4\beta})^N&\sum_{Q\in\Mcal_N}\Big(\prod_{k=1}^N\frac{|\L_N|^{N\widehat Q(k)}}{(N\widehat Q(k))!k^{N\widehat Q(k)}(4\pi\beta k)^{d/2N\widehat Q(k)}}\Big)\le Z_N^{\ssup{\rm sym}}(\beta)\\ &\le \sum_{Q\in\Mcal_N}\Big(\prod_{k=1}^N\frac{|\L_N|^{N\widehat Q(k)}}{(N\widehat Q(k))!k^{N\widehat Q(k)}(4\pi\beta k)^{d/2N\widehat Q(k)}}\Big).
\end{aligned}
$$

Recall from \eqref{functionalS} that
$$
\Scal(Q)=\sum_{k=1}^N\widehat Q(k)\Big[\log\frac{\widehat Q(k)}{\widehat Q^*(k)}-1\Big]\quad \mbox{ for } Q\in\Mcal_N.
$$

Stirling's formula $ N!\approx N^{N+1/2}{\rm e}^ {-N}\sqrt{2\pi}(1+1/12N+O(N^2)) $ gives the following upper and lower  bound

\begin{equation}\label{upperbZ}
\begin{aligned}
Z_N^{\ssup{\rm sym}}(\beta)
&\le \sum_{Q\in\Mcal_N}\Big(\prod_{k=1}^N\frac{{\rm e}^{N\widehat Q(k)}}{\sqrt{2\pi}\rho^{\widehat Q(k)+1/2}\widehat Q(k)^{N\widehat Q(k)}k^{N\widehat Q(k)}(4\pi\beta k)^{d/2N\widehat Q(k)}}\Big)\\
&\le C^{N^\alpha}{\rm e}^{N\log(1+O(N))}\sum_{Q\in\Mcal_N}\exp\Big(-N\Scal(Q)\Big),
\end{aligned}
\end{equation}

and 

\begin{equation}\label{lowerbZ}
\begin{aligned}
Z_N^{\ssup{\rm sym}}(\beta)&\ge (1-{\rm e}^{-dN/4\beta})^N{\rm e}^{-\frac {1}{2}N^\alpha\log N}{\rm e}^{-N\log(1+O(N))}\sum_{Q\in\Mcal_N}\exp\Big(-N\Scal(Q)\Big),
\end{aligned}
\end{equation}
for $ \alpha\in (\frac{1}{2},1) $. The error bounds in \eqref{upperbZ} and \eqref{lowerbZ} due to Stirling's formula follow from the observation that there is a $ \alpha\in (\frac{1}{2},1) $ such $$ \sharp\{k\colon N\widehat Q(k)\ge 1,\sum_{k=1}^Nkr_k=N\} =N^\alpha. $$ 
Recall that $ N\widehat Q(k)=r_k, k=1,\ldots, N, $ the existence of such an $ \alpha $ can be seen from 
$$
\sharp\{k\colon N\widehat Q(k)\ge 1,\sum_{k=1}^Nkr_k=N\}\le \sharp\{k\colon N\widehat Q(k)=1,\sum_{k=1}^Nkr_k=N\}\le c N^{\frac{1}{2}}.
$$

Thus an upper \eqref{upperbZ} and a lower bound \eqref{lowerbZ} is established. To proof Theorem~\ref{THMfree} we need to study the large $ N$-limit for the upper and lower bound, and the continuity of the functional $ \Scal $. 

\begin{prop}\label{prop-proof}

Let $ (\L_N)_{N\ge 1} $ be a sequence of centred cubes with $ \L_N\uparrow\R^d $ 
and $ N/|\L_N|\to\rho $ as $ N\to\infty $. 
\begin{enumerate} 
\item[(a)] Assume additionally that $ \dist(0,\partial\L_N)=cN $ for some $ c>0 $. Then
$$
\liminf_{N\to\infty}\frac{1}{N}\log \sum_{Q\in\Mcal_N}\exp\Big(-N\Scal(Q)\Big)\ge -\inf_{Q\in\Mcal}\Scal(Q).
$$

\item[(b)] 
$$
\limsup_{N\to\infty}\frac{1}{N}\log\sum_{Q\in\Mcal_N}\exp\Big(-N\Scal(Q)\Big)\le -\inf_{Q\in\Mcal}\Scal(Q).
$$
\end{enumerate}
\end{prop}

\begin{proofsect}{Proof}
(a) It suffices to construct, for any $ Q\in\Mcal $, some sequence $ (Q_N)_{N\ge 1} $ with $ Q_N \in \Mcal_N $ which converges strongly to $ Q $ as $ N\to\infty $ such that
\begin{equation}\label{secondstep}
\limsup_{N\to\infty} \Scal(Q_N)\le  S(Q).
\end{equation}
Let $ Q\in\Mcal $ with $ \sum_{k\in\N}Q(k)=1 $ be given. 
Recall that $ \floor{x} $ is the largest integer smaller or equal to $ x\in\R $. For $ \alpha\in(0,1) $ define
\begin{equation}\label{Q-Ndef}
Q_N(k)=\left\{\begin{array}{r@{\,,\,}l}
1-\sum_{k=2}^{\floor{N^\alpha}}\frac{\floor{NQ(k)}}{N} & \mbox{ if } k=1\\[1.5ex]
\frac{\floor{NQ(k)}}{N} & \mbox{ if } k=2,\ldots,\floor{N^\alpha}\\[1.5ex]
0 & \mbox{ otherwise}
\end{array}\right..
\end{equation}
It is clear that $ \sum_{k=1}^N Q_N(k)=1 $, and $ Q_N(l)\ge Q_N(l+1) $ for $ l\ge 2 $ and $ Q_N(1)\ge 1-\sum_{k\ge 2}Q(k)\ge Q(1)\ge Q(2)\ge Q_N(2) $ follows from the fact that $ x\le y $ implies $ \floor{x}\le \floor{y} $ for $ x,y\in\R_+ $. The sequence $ (Q_N)_{N\in\N} $ converges strongly against $ Q $  as $ N\to\infty $ because
$$
\sum_{k=1}^\infty\big|Q(k)-Q_N(k)\big|=\sum_{k=1}^{\floor{N^\alpha}}\big|Q(k)-Q_N(k)\big|+\sum_{k>\floor{N^\alpha}}Q(k)\le N^{\alpha-1}+\sum_{k>\floor{N^\alpha}}Q(k).
$$
To prove \eqref{secondstep} note that the functional $ \Scal $ is lower semi continuous, because $ \Scal $ is the sum of the relative entropy of $ \widehat Q_N $ with respect to $ \widehat Q^* $ restricted to $ \{1,\ldots,\floor{N^\alpha}\} $ and the total mass $ \sum_{k=1}^{\floor{N^\alpha}}\widehat Q_N(k)\le 1 $, and the relative entropy is lower semi continuous. Hence,
$ \liminf_{N\to\infty}\Scal(Q_N)\ge S(Q) $, because of the strong convergence of $ Q_N $ to $ Q $ as $ N\to\infty $. We are left to prove the upper semi continuity \eqref{secondstep} of the functional $ \Scal $ for the sequence $ (Q_N)_{N\ge 1} $. This is seen from
$$
\Scal(Q_N)=\sum_{k=1}^N\Big(\widehat Q_N(k)\log\widehat Q_N(k)-\widehat Q_N(k)(1+\log \widehat Q^*(k))\Big),
$$
where we have to analyse only the entropy term, i.e., the first term on the right hand side. The entropy is defined by $ H(Q)=-\sum_{k=1}^\infty\widehat Q(k)\log\widehat Q(k) $ for $ Q\in\Mcal $. 
\begin{equation}\label{entropyest}
\begin{aligned}
H(Q_N)-H(Q)&=\sum_{k=1}^{\floor{N^\alpha}}\widehat Q_N(k)\log\frac{\widehat Q_N(k)}{\widehat Q(k)}+\sum_{k=1}^{\floor{N^\alpha}}(\widehat Q_N(k)-\widehat Q(k))\log\widehat Q(k)\\
&\quad+\sum_{k=\floor{N^\alpha}+1}^\infty \widehat Q(k)\log\widehat Q(k)
\end{aligned}
\end{equation}
In the first two terms on the right hand side of \eqref{entropyest} it is sufficient to consider only sum indices $ k\in\{1,\ldots, \floor{N^\alpha}\} $ with $ \widehat Q(k)\ge 1/N $. The terms for $ k=1 $ in the two sums on the right hand side vanish due to the strong convergence.
The remaining first term on the right hand side of \eqref{entropyest} is bounded by
$$
\sum_{k=2}^{\floor{N^\alpha}}|\widehat Q_N(k)-\widehat Q(k)|\le 2N^{\alpha-1}\to 0\mbox{ as } N\to\infty,
$$
the second term is bounded
$$
\Big|\sum_{\heap{k=1}{\widehat Q(k)>1/N}}^{\floor{N^\alpha}}(\widehat Q(k)_N-\widehat Q(k))\log\widehat Q(k)\Big|\le 2N^{\alpha-1}\log N\to 0\mbox{ as }N\to\infty.
$$
The third term clearly vanishes in the limit $ N\to\infty $.
The latter estimations together with \eqref{entropyest} finish the proof for the upper semi continuity, i.e.,
$ \limsup_{N\to\infty}\Scal(Q_N)\le S(Q) $, and therefore \eqref{secondstep} is proven. This together with the lower bound in \eqref{lowerbZ}, i.e., all the error terms vanish in the limit, proves  (a).

\noindent (b) For the upper bound we need to estimate the cardinality of the set $ \Mcal_N $.
The set $ \Mcal_N $ consists of all antitone mappings $ Q\colon\{1,\ldots,N\}\to\{0,1/N,2/N,\ldots,1\} $ with $ \sum_{l=1}^NQ(l)=1 $ which are elements in the image $ Q_N(\Pcal_N) $ of the discrete empirical shape measure $ Q_N $ \eqref{discreteshapemeasure}. For $ \l\in\Pcal_N $ there are at most $ N^{\frac{1}{2}} $ indices $ k $ such that the occupation number for a component/cycle of length $k $ is nonzero, i.e., $ r_k\ge 1 $. These indices are called occupied. One maximises the number of occupied indices if one starts with index $1$ and fills each component/cycle only ones. Then, at index $ \ceil{N^{1/2}} $ we already exceed the total number $ N $, because $ \sum_{k=1}^Nkr_k =N $. Here $ \ceil{x} $ is the smallest integer greater or equal to $ x\in\R_+ $. Hence, there is a $ \alpha\in(\frac{1}{2},1) $ such that
$$
\sharp\{k\in\{1,\ldots,N\}\colon r_k(\l)\ge 1\}=\sharp\{k\in\{1,\ldots,N\}\colon N\widehat Q^\l(k)\ge 1\}\le N^{\alpha}\quad \mbox{ for } \lambda\in\Pcal_N.
$$
We choose $ \ceil{N^{\alpha}} $ indices out from $ N $ and we denote this set by $ M $. $ M $ has cardinality $ \ceil{N^\alpha} $. The cardinality of the set $ \{f\colon M\to\{1/N,2/N,\ldots,1\}\colon f(a)\ge f(b) \mbox{ for } a\le b, a,b\in M\} $ is then given by
\begin{equation}\label{e1}
\begin{aligned}
\sharp \Big\{f\colon M\to\{1/N,2/N,\ldots,1\}\colon f(a)\ge f(b) \mbox{ for } a\le b, a,b\in M\Big\}&=\frac{N(N+1)\cdots(N+\ceil{N^\alpha}-1)}{\ceil{N^\alpha}!}\\
& =\frac{(N+\ceil{N^\alpha}-1)!}{(N-1)!\ceil{N^\alpha}!}.
\end{aligned}
\end{equation}
Recall that any element in the image $ Q_N(\Pcal_N) $ is antitone and a probability measure, whereas the latter condition is relaxed in \eqref{e1}. We prove \eqref{e1}. Let the set $ M=\{m_1,\ldots,m_{\ceil{N^\alpha}}\} $ be given with the natural order. Then the mapping $ \Phi\colon b_1b_2\cdots b_{\ceil{N^\alpha}}\to\{b_1+(\ceil{N^\alpha}-1), b_2+(\ceil{N^\alpha}-2),\ldots, b_{\ceil{N^\alpha}}\} $ is a bijection between the set of all antitone mappings $ f\colon M\to\{1/N,2/N,\ldots,1\} $ and the set of all $ \ceil{N^\alpha} $-subsets of $ \{1,2,\ldots, \ceil{N^\alpha}+N-1\} $. Here we used the word representation of the mapping $ \Phi $. $ N(N+1)\cdots(N+\ceil{N^\alpha}-1) $ counts the number of sortings of $ \ceil{N^\alpha} $ objects into $ N $ linearly ordered boxes.

The occupation numbers $ r_k=0 $ for all the remaining indices $ k\in\{1,\ldots,N\}\setminus M $. Thus $ Q^\l_N(l)=0 $ for all $ l\in \{1,\ldots,N\}\setminus M $ with $ l>\max M $. Otherwise $ Q^\l_N(l) $ is constant determined by the indices of the set $ M $, i.e., $ Q^\l_N(l)=1/N\sum_{k\in M, k>l}r_k(\l) $ for $ \{1,\ldots,N\}\setminus M $ and $ l<\max M $. Hence, to estimate the cardinality of the image set $ \Mcal_N $ we have $ \binom{N}{\ceil{N^\alpha}} $ ways to choose a subset $ M\subset\{1,\ldots,N\} $ of cardinality $ \sharp M=\ceil{N^\alpha} $. The remaining choices for the image elements are then given by \eqref{e1} because all the indices in the complement $ \{1,\ldots,N\}\setminus M $ are determined by the indices in $ M $.
Hence, using Stirling's formula again, we get
\begin{equation}\label{estM}
\sharp\Mcal_N \le \binom{N}{\ceil{N^\alpha}}\frac{(N+\ceil{N^\alpha}-1)!}{(N-1)!\ceil{N^\alpha}!}\le C N^{2\ceil{N^\alpha}(1-\alpha)}\Big(1-\frac{\ceil{N^\alpha}}{N}\Big)^{2\ceil{N^\alpha}}
\end{equation}
for some $ C>0 $.

Let $ Q_N\in\Mcal_N $ be minimal for $ \Scal $ on the set $ \Mcal_N $ and let $ Q\in\Mcal $ denote the limit $ \lim_{N\to\infty} Q_N$. From the upper bound \eqref{upperbZ} we obtain
\begin{equation}\label{lastupperbZ}
\begin{aligned}
Z_N^{\ssup{\rm sym}}(\beta)
&\le C N^{2\ceil{N^\alpha}(1-\alpha)}\Big(1-\frac{\ceil{N^\alpha}}{N}\Big)^{2\ceil{N^\alpha}}\exp\Big(-N\Scal(Q_N)\Big).
\end{aligned}
\end{equation}
By the continuity of $ Q\mapsto \Scal(Q) $, i.e., by the lower semi continuity of $ \Scal $, the assertion (b) of Proposition~\ref{prop-proof} follows.

\qed
\end{proofsect}

\subsection{Proof of Theorem~\ref{THM-variational}}\label{proof-sec-variational}
We analyse the variational formula in \eqref{variationalfreeenergy} for the functional $ \Scal $ and prove Theorem~\ref{THM-variational}.
Note that any montonuously non increasing probability measure $ Q\in\Mcal $ is in one to one correspondence to $ \widehat Q\in \{\widehat Q\colon\N\to [0,1]\colon \sum_{k=1}^\infty \widehat Q(k)\le 1\} $. If $ Q,P\in\Mcal $ the equality $ \widehat Q=\widehat P $ implies $ Q(1)-Q(n)=P(1)-P(n) $ for $ n\in\N $, which implies $ Q(1)=P(1) $ and hence $ P=Q $. Furthermore note that the functional $ \Scal $ is convex and lower semi continuous and can be written as 
\begin{equation}\label{relentropy}
S(Q)=qH(\widehat P|\widehat P^*)+q\log\frac{q}{q^*}-q, 
\end{equation} where $ q=\sum_{k=1}^\infty\widehat Q(k) , q*=\sum_{k=1}^\infty\widehat Q^*(k) $ and $ H $ the relative entropy of the two probability measures $ \widehat P=q^{-1}\widehat Q $ and $ \widehat P^*=(q^*)^{-1} \widehat Q^* $. 
A glance at the first term in \eqref{relentropy} would imply that $ \widehat P=\widehat P^* $, because the relative entropy is zero if and only if $ \widehat P=\widehat P^* $.

For $ d=1,2 $, and any density $ \rho<\infty $, there is a unique solution $ \alpha $ of \eqref{root}. Hence, $ Q\in\Mcal $ with 
$$
\widehat Q(k)=\frac{{\rm e}^{-\alpha k}}{\rho(4\pi\beta)^{d/2}k^{1+\frac{d}{2}}}\quad\mbox{ for } k\in\N,
$$ 
is the unique solution of the variational formula \eqref{variationalfreeenergy} with $ \sum_{k=1}^\infty Q(k)=1 $. The free energy is given by $ \rho/\beta\Scal(Q) $, and \eqref{freeed12} is proved.

We turn to the case $ d\ge 3 $ in the following. 
We fix $ \rho<\rho_c $. Minimising only the relative entropy in \eqref{relentropy} is not sufficient because $ Q=Q^* $ would imply that
$$
\sum_{k=1}^\infty Q(k)=\sum_{k=1}^\infty k\widehat Q(k)=\frac{1}{\rho(4\pi\beta)^{d/2}}\zeta(d/2)=\frac{1}{\rho(4\pi\beta)^{d/2}}g_{d/2}(0) >1.
$$
This follows from the fixed density in \eqref{root}, i.e., $ \rho=\frac{1}{(4\pi\beta)^{d/2}}g_{d/2}(\alpha) $ and $ g_{d/2}(\alpha)<g_{d/2}(0) $ (see e.g. \cite{Gram25}). Hence the unique solution for the variational problem in \eqref{variationalfreeenergy} is given by $ Q\in\Mcal $ with
\begin{equation}
\widehat Q(k)=\frac{{\rm e}^{-\alpha k}}{\rho(4\pi\beta)^{d/2}k^{1+\frac{d}{2}}}\quad\mbox{ for } k\in\N,
\end{equation}
and the free energy is given as
$$
\frac{\rho}{\beta}\sum_{k=1}^\infty\frac{{\rm e}^{-\alpha k}}{k^{1+d/2}\rho(4\pi\beta)^{d/2}}(-\alpha k-1)=-\frac{1}{\beta(4\pi\beta)^{d/2}}g_{\frac{d+2}{2}}(\alpha)-\frac{1}{\beta}\alpha\rho
$$ for $ \rho<\rho_c $.

For the case $ \rho>\rho_c $ define the sequence $ (Q_n)_{n\in\N} $ by
\begin{equation}
\widehat Q_n(k)=\left\{\begin{array}{r@{\;,\;}l}
\widehat Q^*(k) & \mbox{ for } k\not= n\\
\widehat Q^*(k)+\frac{(\rho-\rho_c)}{n\rho} & \mbox{ for } k=n
                       \end{array}\right.
\end{equation}
Recall that $ Q_n(l)=Q_n(1)-\sum_{k=1}^{l-1}\widehat Q_n(k) $ for any $ l\in\N $. Then $ Q_n\in\Mcal, n\in\N $, because $ Q_n(k)\ge Q_n(k+1), k\in\N, $  and 
$$ 
\sum_{k=1}^\infty Q_n(k)= \sum_{k=1}^\infty k\widehat Q_n(k)=\frac{1}{\rho(4\pi\beta)^{d/2}}\zeta\Big(\frac{d}{2}\Big)+\frac{(\rho-\rho_c)}{\rho}=1.
$$

We have $ Q_n\to Q^* $ strongly as $ n\to\infty $ and for any $ n\in\N $ the evaluation for the functional is given by
$$
S(Q_n)=S(Q^*)-\frac{(\rho-\rho_c)}{n\rho}+\Big(\widehat Q^*(n)+\frac{(\rho-\rho_c)}{n\rho}\Big)\log\Big(\frac{\widehat Q^*(n)+\frac{(\rho-\rho_c)}{n\rho}}{\widehat Q^*(n)}\Big).
$$
Hence, $ \lim_{n\to\infty}S(Q_n)=S(Q^*) $. The infimum is therefore not attained within the set $ \Mcal $ for $ \rho>\rho_c $, only minimising sequences do exist.

\subsection{Proof of Theorem~\ref{THMLDP-LN}}\label{proof-secLDP-LN}
The following is a reformulation of Theorem~\ref{THMLDP-LN}, making explicit what a large deviations principle is.
Let $ \L_N\subset\R^d $ with $\dist(0,\partial\L_N)=cN $ for some $ c>0 $ and $ N/|\L_N|\to\rho\in (0,\infty) $ and $ \L_N\uparrow\R^d $ as $ N\to\infty $, then Theorem~\ref{THMLDP-LN} is equivalent to (i)-(iii):

\noindent {\bf (i) Lower bound of the large deviations principle}

For any open set $ G\subset\Mcal_1(\Omega) $,
\begin{equation}\label{lowerTHM}
\liminf_{N\to\infty}\frac{1}{N}\log\P_N^{\ssup{\rm sym}}(L_N\in G)\ge -\inf_{\mu\in G} I^{\ssup{\rm sym}}(\mu).
\end{equation}

\noindent {\bf (ii) Upper bound of the large deviations principle}

For any compact set $ K\subset\Mcal_1(\Omega) $,
\begin{equation}\label{upperTHM}
\limsup_{N\to\infty}\frac{1}{N}\log\P_N^{\ssup{\rm sym}}(L_N\in K)\le -\inf_{\mu\in K} I^{\ssup{\rm sym}}(\mu).
\end{equation}

\noindent {\bf (iii) Exponential tightness}

The sequence of distributions of $ L_N $ under $ \P_N^{\ssup{\rm sym}} $ is exponentially tight.

\noindent The proof of (iii) is given in Lemma~\ref{exptightnesssym} in Subsection~\ref{sec-exptight}. The proof of the lower bound \eqref{lowerTHM} and upper bound \eqref{upperTHM} follows below. We are going to estimate the probabilities directly.  The key idea is to split the two random mechanism. First, rewrite the sum over permutations in a sum over discrete empirical shape measures for integer partitions. The asymptotic of the discrete empirical shape measure is governed by the functional $ \Scal $ \eqref{functionalS}. Any discrete empirical shape measure $ Q\in\Mcal_N $ determines a probability measure $ \P^{\ssup{Q}}_N $ on $ \Omega_\beta^N $ (respectively on $ \Omega $), which represents the concatenations of Brownian motions with time horizon $ [0,\beta] $ to Brownian bridges with time horizon an integer multiple of $ [0,\beta] $ depending on the length of the cycle present for the discrete empirical shape measure $ Q $. There is a large deviations principle for $ L_N $ under $ \P_N^{\ssup{Q}} $ with rate function $ I^{\ssup{Q}} $ \eqref{defrateQ}.

\begin{proofsect}{Proof of the lower bound \eqref{lowerTHM} of Theorem~\ref{THMLDP-LN}}
Let $ G\subset\Mcal_1(\Omega) $ open with respect to the weak topology for probability measure. As $ \Omega $ is Polish, the L\`{e}vy metric is compatible with this topology (\cite{DS01}). Replace the sum over permutations in the definition \eqref{symPdef1} of the symmetrised measure $ \P_N^{\ssup{\rm sym}} $ by a sum over partitions in $ \Pcal_N $. As in the proof of Theorem~\ref{THMfree} in Subsection~\ref{proof-sec-free} we regroup for a given partition $ \l\in\Pcal_N $ the product of the Brownian bridges. In the following we write $ r_k $ for $ r_k(\l) $. Recall the definition of the measure $ \mu_k=\mu_{0,0}^{\beta,N} $ in \eqref{defmuk}, and define for any partition $ \l\in\Pcal_N $ the measure
$ \bigotimes_{k=1}^N\P_k^{\otimes r_k} \in\Mcal_1(\Omega_\beta^N)$,
where  $ \P_k $ is the normalised version of the measure $ \mu_k $, i.e.,
\begin{equation}
\P_k=\mu_k/\mu_k(\Omega_\beta^N=\mu_k/p_{k\beta,\L_N}(0,0).
\end{equation}
Hence,
$$
\begin{aligned}
\P_{N}^{\ssup{\rm sym}} &=\frac{1}{Z_N^{\ssup{\rm sym}}(\beta)}\sum_{\l\in\Pcal_N}\Big(\prod_{k=1}^N\frac{|\L_N|^{r_k}g_{k\beta}(0,0)^{r_k}}{r_k!k^{r_k}}\Big)\bigotimes_{k=1}^N\P_k^{\otimes r_k},
\end{aligned}
$$
where the volume factor in the numerator comes from the volume factor in the definition of the measure $ \mu_k $ in \eqref{defmuk}. Thus,

$$
\begin{aligned}
\P_{N}^{\ssup{\rm sym}} &=\frac{1}{Z_N^{\ssup{\rm sym}}(\beta)}\sum_{Q\in\Mcal_N}\Big(\prod_{k=1}^N  \frac{|\L_N|^{N\widehat Q(k)}g_{k\beta}(0,0)^{N\widehat Q(k)}}{ (N\widehat Q(k))!k^{N\widehat Q(k)}}\Big)\P^{\ssup{Q}}_{N},
\end{aligned}
$$
where
\begin{equation}\label{def-P-Q}
\P^{\ssup{Q}}_{N}=\bigotimes_{k=1}^N\P_k^{\otimes N\widehat Q(k)}\in\Mcal_1(\Omega_\beta^{\otimes N}),\quad Q\in\Mcal_N.
\end{equation}

We now prove the lower bound for the large deviations principle in Theorem~\ref{THMLDP-LN}. $ \P_N^{\ssup{Q}}\circ L_N^{-1} $ is then a distribution in $ \Mcal_1(\Omega) $, i.e., we conceive $ \P_N^{\ssup{Q}} $ as measure on $ \Omega $.

Stirling's formula gives with the lower bound \eqref{lowerbZ} the following lower bound
\begin{equation}\label{problower1}
\begin{aligned}
\P_{N}^{\ssup{\rm sym}}&(L_N\in G)\ge(1-{\rm e}^{-dN/4\beta})^N{\rm e}^{-\frac {1}{2}N^\alpha\log N}{\rm e}^{-N\log(1+O(N))}\\ &\times\frac{1}{Z_N^{\ssup{\rm sym}}(\beta)} \sum_{Q\in\Mcal_N}\Big(\prod_{k=1}^N\frac{{\rm e}^{N\widehat Q(k)}}{\rho^{\widehat Q(k)}\widehat Q(k)^{N\widehat Q(k)}k^{N\widehat Q(k)}(4\pi\beta k)^{d/2N\widehat Q(k)}}\Big)  \bigotimes_{k=1}^N \P_k^{\otimes N\widehat Q(k)}(L_N\in G)\\=
&(1-{\rm e}^{-dN/4\beta})^N{\rm e}^{-\frac {1}{2}N^\alpha\log N}{\rm e}^{-N\log(1+O(N))} \frac{1}{Z_N^{\ssup{\rm sym}}(\beta)}\sum_{Q\in\Mcal_N}\exp\Big(-N\Scal(Q)\Big)\P_N^{\ssup{Q}}(L_N\in G),
\end{aligned}
\end{equation}
where $ \alpha\in (\frac{1}{2},1) $, and where $ \dist(0,\partial\L_N)=cN $ for some $ c>0 $ (see Proof of Theorem~\ref{THMfree}).

We prove the large deviations lower bound for a sequence of boxes $ \L_N\subset\R^d $ with $ \L_N\uparrow\R^d $ as $ N \to\infty $ and $ \dist(0,\partial\L_N)=cN $ for some $ c>0 $.

\begin{prop}\label{propproof2} Let $ (\L_N)_{N\ge 1} $ be a sequence of centred cubes with $ \L_N\uparrow\R^d $ 
and $ N/|\L_N|\to\rho $ as $ N\to\infty $ and $ \dist(0,\partial)=cN $ for some $ c>0 $.
Then
$$
\liminf_{N\to\infty}\frac{1}{N}\log \Big(\sum_{Q\in\Mcal_N}\exp\Big(-N\Scal(Q)\Big)\P_N^{\ssup{Q}}(L_N\in G)\ge -\inf_{\mu\in A} I^{\ssup{\rm sym}}(\mu).
$$

\end{prop}

\begin{proofsect}{Proof of Proposition~\ref{propproof2}}
It suffices to construct, for any $ \mu\in G $ and $ Q\in\Mcal $, some sequence $ (Q_N)_{N\ge 1} $ with $ Q_N\in\Mcal_N $ and with $ Q_N\to Q $ as $ N\to\infty $ strongly such that
\begin{equation}\label{step3}
\limsup_{N\to\infty}\Big(\Scal(Q_N)-\frac{1}{N}\log\P^{\ssup{Q_N}}(L_N\in G)\Big)\le \Scal(Q) +I^{\ssup{Q}}(\mu).
\end{equation}
The first part of \eqref{step3}, i.e., $ \limsup_{N\to\infty}\Scal(Q_N)\le \Scal(Q) $  is shown in the proof of Proposition~\ref{prop-proof} in \eqref{secondstep}. 
By $ \d $  we denote the L\'{e}vy metric on $\Mcal_1(\Omega)$, which generates the weak topology; see \eqref{Levymetric} below. 

\noindent By $\dist(\mu,G)=\inf_{\nu\in G}\d(\mu,\nu)$ we denote the distance to a set $G\subset \Mcal_1(\Ccal) $. Recall the L\'{e}vy metric $ \d $ on the Polish space $ \Mcal_1(\Omega) $ \cite{DS01}, defined for any two probability measures $ \mu,\nu\in\Mcal_1(\Omega) $ as
\begin{equation}\label{Levymetric}
\d(\mu,\nu)=\inf\{\delta >0\colon\mu(\Gamma)\le \nu(\Gamma^{\delta})+\delta\,\mbox{ and }\, \nu(\Gamma)\le\mu(\Gamma^{\delta})+\delta\,\mbox{ for all }\, \Gamma=\overline{\Gamma}\subset\Omega\}, 
\end{equation} 
where $ F^{\delta}=\{\mu\in\Mcal_1(\Omega)\colon\dist(\mu,F)\leq\delta\} $ is the closed $\delta$-neighbourhood of $ F $. 
In the following let $ \delta_N $ such that $ \delta_N\downarrow 0 $ as $ N\to\infty $. Let $ G_{\delta_N} $ be a $ \delta_N$-neighbourhood of $ \mu $ in $ G $, i.e. $ G_{\delta_N}=\{\nu\in\Mcal_1(\Omega)\colon \dist(\nu,\mu)\le\delta_N\}$.

Let $ Q_N \in\Mcal_N $ defined as in \eqref{Q-Ndef}. We are going to use the G\"artner-Ellis Theorem to deduce that
$$
\liminf_{N\to\infty}\frac{1}{N}\log \P^{\ssup{Q_N}}(L_N\in G_{\delta_N})\ge I^{\ssup{Q}}(\mu). 
$$

For doing this, we first introduce for any $ F\in\Ccal_{\rm b}(\Omega) $ the logarithmic moment generating function
$$
\begin{aligned}
\L_N(F)&:=\log\E_{\P_N^{\ssup{Q_N}}}\Big({\rm e}^{N\langle F,L_N\rangle}\Big)=\log\Big(\prod_{k=1}^N\E_{0,0}^{k\beta}\Big({\rm e}^{F(B)}\Big)^{N\widehat Q_N(k)}\Big)\\
&=N\sum_{k=1}^N\widehat Q_N(k)\log\E_{0,0}^{k\beta}\Big({\rm e}^{F(B)}\Big),
\end{aligned}
$$
where $ \E_{\P_N^{\ssup{Q_N}}} $ denotes the expectation with respect to the probability measure $ \P_N^{\ssup{Q}} $, and $ \E_{0,0}^{k\beta} $ the expectation with respect to the probability measure $ \P_k $. Here $ \P_k $ is supported on the path set $ A_k $ \eqref{defAK}.
From the strong convergence of $ Q_N $ to $ Q $ as $ N \to\infty $ it is easily seen that
$ \L(F):=\lim_{N\to\infty}\frac{1}{N}\L_N(F) $ exists, and 
$$
\L(F)=\sum_{k\in\N}\widehat Q(k)\log\E_{0,0}^{k\beta}\Big({\rm e}^{F(B)}\Big).
$$
Since it is easily seen that $ \L $ is lower semi continuous and G\^ateaux differentiable, and by the exponential tightness of the family $ \P^{\ssup{Q_N}}_N\circ L_N^{-1} $ (see Lemma~\ref{exptightnessQ} in Subsection~\ref{sec-exptight}), \cite[4.5.27]{DZ98} implies that
$$
\liminf_{N\to\infty}\frac{1}{N}\log \P_N^{\ssup{Q_N}}(L_N\in G_{\delta_N})\ge - I^{\ssup{Q}}(\mu).
$$
This shows the second half of  \eqref{step3}and finishes the proof.
\qed
\end{proofsect}

\begin{proofsect}{Proof of the upper  bound \eqref{upperTHM} of Theorem~\ref{THMLDP-LN}}
From the exponential tightness (iii) in Lemma~\ref{exptightnessQ} it follows that there is a sequence of compact sets $ 
K_L\subset\Mcal_s(\Omega) $ such that
$$
\lim_{L\to\infty}\limsup_{N\to\infty}\frac{1}{N}\log\Big(\sup_{Q\in\Mcal_N}\P_N^{\ssup{Q}}(L_N\in K_L^{\rm c})\Big)=-\infty.
$$
Hence, for the proof of the upper bound of the large deviations principle it suffices to show \eqref{upperTHM} for compact sets.
Let $ K\subset\Mcal_1(\Omega) $ compact. Then \eqref{upperbZ} and \eqref{def-P-Q} imply
\begin{equation}
\begin{aligned}
\P_N^{\ssup{\rm sym}}(L_N\in K)\le C^{N^\alpha}{\rm e}^{N\log(1+O(N))}\sum_{Q\in\Mcal_N}\exp\Big(-N\Scal(Q)\Big)\P_N^{\ssup{Q}}(L_N\in K)
\end{aligned}
\end{equation}

We consider the logarithmic moment generating function of the distribution of $ L_N $ under $ \P_N^{\ssup{Q}} $,
$$
\begin{aligned}
\L_N(F):=\log\E_{\P_N^{\ssup{Q_N}}}\Big({\rm e}^{N\langle F,L_N\rangle}\Big)=N\sum_{k=1}^N\widehat Q_N(k)\log\E_{0,0}^{k\beta}\Big({\rm e}^{F(B)}\Big)\quad , F\in\Ccal_{\rm b}(\Omega). 
\end{aligned}
$$

Let now $ Q_N\in\Mcal_N $ be maximal for 
$$ Q\mapsto{\rm e}^{-N\Scal(Q)}\P_N^{\ssup{\rm sym}}(L_N\in K). 
$$ 
The sequence $ (Q_N)_{N\ge 1} $ has a limit $ Q\in\Mcal $ as $ N\to\infty $. Certainly, the limit
$ \L(F):=\lim_{N\to\infty}\frac{1}{N}\L_N(F) $ exists, and is lower semi continuous and G\^ateaux differentiable. Observe that $ I^{\ssup{Q}} $ is the Fenchel-Legendre transform of $ \L $. Now the G\"artner-Ellis Theorem yields that
$$
\limsup_{N\to\infty}\frac{1}{N}\P_N^{\ssup{Q_N}}(L_N\in K)\le -\inf_{\mu\in K} I^{\ssup {Q}}(\mu).
$$
The cardinality of $ \Mcal_N $ is bounded by \eqref{estM}, and by the continuity of $ Q\mapsto\Scal(Q) $, the assertion follows.

\qed
\end{proofsect}

\subsection{Exponential tightness}\label{sec-exptight}
\noindent In this subsection, we prove the necessary exponential tightness assertions for the sequence of the empirical path measures under the symmetrised measures,  $ \P_{N}^{\ssup{\rm sym}} $, and  under the mixed product measures, $  \P^{\ssup{Q}}_N, $ for any $ Q\in\Mcal_N $. The proof of the latter exponential tightness is a variant of the standard proof for laws of empirical measures. Here, the main ingredient is the product structure of the probability measure. 
We equip the path space $ \Omega $ with the following complete metric
\begin{equation}
\d(\omega,\xi)=\sum_{k=1}^\infty2^{-k}\d_k(\omega,\xi),\quad\omega,\xi\in\Omega,
\end{equation}
where 
$$ 
\d_k(\omega,\xi)=\big\|(\omega-\xi)|_{[0,k\beta]}\big\|_{\infty}\wedge 1,\quad k\in\N. 
$$

Then $ \d $ is a complete metric on $ \Omega $, which generates the uniform convergence on compacts, and $ (\Omega,\d) $ is a Polish space.

\begin{lemma}\label{exptightnesssym}
The family of distributions of the empirical path measures $ L_N $ under the symmetrised measure $ \P_N^{\ssup{\rm sym}} $ is exponentially tight.
\end{lemma}

\begin{proofsect}{Proof}
The proof is in the spirit of the proof of \cite[Lemma~6.2.6]{DZ98}. Choose  $ \delta_l , l\in\N $, so small that for all $ k\in\N $,
\begin{equation}
 \P_{0,0}^{k\beta}\Big(\sup_{\heap{s,t\in[0,k\beta],}{|s-t|\le\delta_l}}|B_s-B_t|>\frac{1}{l}\Big)\le {\rm e}^{-l^2}.
\end{equation}
Consider 

$$
A_l=\Big\{\omega\in\Omega\colon \omega(0)=0,\forall k\in\N\colon \sup_{\heap{s,t\in[0,k\beta],}{|s-t|\le\delta_l}}|\omega_s-\omega_t|\le\frac{1}{l}\Big\}.
$$

According to Arzel\`a-Ascoli's Theorem, $ A_l $ is relative compact in $ \Omega $, because the set $ \{\omega(0)\colon\omega\in A_l\} $ of the starting points is bounded,  and for any $ k\in\N $ we have
$$
\lim_{l\to\infty}\sup_{\omega\in A_l}\sup_{\heap{s,t\in[0,k\beta],}{|s-t|\le\delta_l}}|\omega_s-\omega_t|= 0.
$$
Put $ M_l:=\{\mu\in\Mcal_1(\Omega)\colon\mu(\overline{A}_l^{\rm c})\le \frac{1}{l}\} $ and note that $ M_l $ is closed by Portmanteau's theorem. Let $ L\in\N $ be given and consider $ K_L:=\bigcap_{l=L}^{\infty} M_l $. It is easy to see that $K_L$ is tight, hence $\overline K_L$ is compact by Prohorov's theorem. We shall show that $ \P_{N}^{\ssup{\rm sym}}(L_N\in \overline K_L^{\rm c})\le {\rm e}^{-LN} $ for any $N\in\N$, which implies the assertion. Observe that 
$$
\{L_N\in M_l^{\rm c}\}\subset\Big\{\sharp\{i\in\{1,\ldots,N\}\colon B^{\ssup{i}}\otimes_\beta P^{B^{\ssup{i}}_\beta}\in A_l^{\rm c}\}>\frac{N}{l}\Big\}.
$$
Then, using the upper bounds in \eqref{upperbZ} and \eqref{lastupperbZ}, we get

\begin{equation}
\begin{aligned}
\P_N^{\ssup{\rm sym}}&(L_N\in M_l^{\rm c})\le \P_N^{\ssup{\rm sym}}\Big(\sharp\Big\{\sup_{\heap{s,t\in[0,\beta],}{|s-t|\le\delta_l}}|B^{\ssup{i}}_s-B^{\ssup{i}}_t|>\frac{1}{l}\Big\}\ge \frac{N}{l}\Big)\\
&\le \sum_{\heap{I\subset\{1,\ldots,N\},}{|I|\ge N/l}}\sum_{Q\in\Mcal_N}\frac{{\rm e}^{-N\Scal(Q)}}{Z_N^{\ssup{\rm sym}}(\beta)}\bigotimes_{k=1}^N\P_k\Big(\forall i\in I\colon B_0^{\ssup{i}}=0,\sup_{\heap{s,t\in[0,\beta],}{|s-t|\le\delta_l}}|B^{\ssup{i}}_s-B^{\ssup{i}}_t|>\frac{1}{l}\Big)\\
& \le \sum_{\heap{I\subset\{1,\ldots,N\},}{|I|\ge N/l}}\sum_{Q\in\Mcal_N}\frac{{\rm e}^{-N\Scal(Q)}}{Z_N^{\ssup{\rm sym}}(\beta)}\bigotimes_{\heap{k\in\{1,\ldots,N\},}{kN\widehat Q(k)\ge |I|}}\P_{0,0}^{k\beta}\Big(\sup_{\heap{s,t\in[0,k\beta],}{|s-t|\le\delta_l}}|B_s-B_t|>\frac{1}{l}\Big)^{N\widehat Q(k)}\\
&\le \sum_{\heap{I\subset\{1,\ldots,N\},}{|I|\ge N/l}}\sum_{Q\in\Mcal_N}\frac{{\rm e}^{-N\Scal(Q)}}{Z_N^{\ssup{\rm sym}}(\beta)}\prod_{\heap{k\in\{1,\ldots,N\},}{kN\widehat Q(k)\ge |I|}}{\rm e}^{-l^2kN\widehat Q(k)}\\
&\le C N^{2\ceil{N^\alpha}(1-\alpha)}\Big(1-\frac{\ceil{N^\alpha}}{N}\Big)^{2\ceil{N^\alpha}}\sum_{\heap{I\subset\{1,\ldots,N\},}{|I|\ge N/l}}{\rm e}^{-lN}\le CN^{O(N)}2^N {\rm e}^{-Nl}.
\end{aligned}
\end{equation}
Hence,
$$
\P_N^{\ssup{\rm sym}}(L_N\in K_L^{\rm c})\le\sum_{l=L}^{\infty}\P_{N}^{\ssup{\rm sym}}(L_N\in M^{\rm c}_l)\le CN^{O(N)}2^N\sum_{l=L}^{\infty}{\rm e}^{-lN}\le CN^{O(N)}2^{N+1}{\rm e}^{-NL}\le {\rm e}^{-NL/R}
$$ for all large $ N $ if $ R $ is chosen with
$$ 
CN^{O(N)}2^{N+1}\le{\rm e}^{NL(1-1/R)}.
$$ This ends the proof.

\qed
\end{proofsect}
Now we prove the exponential tightness of the empirical path measures $L_N$ under the probability measures $  \P^{\ssup{Q}}_N, Q\in\Mcal_N $, introduced in \eqref{def-P-Q}.
\begin{lemma}\label{exptightnessQ}
Let $ (Q_N)_{N\ge 1} $ be a sequence of measures $ Q_N\in\Mcal_N $ such that $ Q_N\to Q\in\Mcal $ strongly as $ N\to\infty $. Then the family of distributions of the empirical path measures $ L_N $ under the measures $ \P_N^{\ssup{Q_N}} $ is exponentially tight.
\end{lemma}

\begin{proofsect}{Proof}
As in the proof of Lemma~\ref{exptightnesssym} there exists $ l\in\N $ and compact sets $ \Gamma_l\subset\Omega $ such that for all $ k\in\N $
\begin{equation}
\P_k(B\otimes_{k\beta}P^0\in\Gamma_l^{\rm c})\le {\rm e}^{-2l^2}({\rm e}^l-1).
\end{equation}

The set $ M_l=\{\nu\in\Mcal_1(\Omega)\colon \nu(\Gamma_l^{\rm c})\le 1/l\} $ is closed by  Portmanteau's theorem. For $ L\in\N $ define $ K_L=\bigcap_{l=L}^{\infty}M_l $. By Prohorov's theorem, each $ K_L $ is a relative compact subset of $ \Mcal_1(\Omega) $. Then Chebycheff's inequality gives that for any $ N\in\N $, and any $ Q_N\in\Mcal_N $  
\begin{equation}
\begin{aligned}
\P_N^{\ssup{Q_N}}(L_N\notin M_l)&=P_N^{\ssup{Q_N}}(L_N(\Gamma_l^{\rm c})>\frac{1}{l})\le \E_{\P_N^{\ssup{Q_N}}}\Big({\rm e}^{2Nl^2(L_N(\Gamma_l^{\rm c})-1/l)}\Big)\\
&={\rm e}^{-2Nl}\E_{\P_N^{\ssup{Q_N}}}\Big(\exp\Big(2l^2\sum_{i=1}^N\1\{B^{\ssup{i}}\otimes_\beta P^{B_\beta^{\ssup{i}}}\in\Gamma_l^{\rm c}\}\Big)\Big)\\
&={\rm e}^{-2Nl}\prod_{k=1}^N\E_{0,0}^{k\beta}\Big(\exp\Big(2l^2\1\{B^{\ssup{i}}\otimes_{k\beta} P^{B_\beta}\in\Gamma_l^{\rm c}\}\Big)\Big)^{N\widehat Q_N(k)}\\
&={\rm e}^{-2Nl}\prod_{k=1}^N\Big(\P_k(B\otimes_{k\beta}P^{B_{k\beta}}\in\Gamma_l)+{\rm e}^{2l^2}\P_k(B\otimes_{k\beta}P^{B_{k\beta}}\in\Gamma_l^{\rm c})\Big)^{N\widehat Q_N(k)}\\
&\le {\rm e}^{-2Nl}\prod_{k=1}^N({\rm e}^l)^{N\widehat Q_N(k)}\le {\rm e}^{-Nl}.
\end{aligned}
\end{equation}
Therefore,
\begin{equation}
\P_N^{\ssup{Q_N}}(L_N\notin K_L)\le\sum_{l=L}^\infty\P_N^{\ssup{Q_N}}(L_N\notin M_l)\le\sum_{l=L}^\infty{\rm e}^{-Nl}\le 2{\rm e}^{-NL},
\end{equation}
which implies the exponential tightness.
\qed
\end{proofsect}

\end{proofsect}
\section{Appendix}\label{sec-app}
\subsection{Large deviations}\label{sec-ldp}
\noindent For the convenience of the reader, we repeat the notion of a large deviations principle and of the most important facts that are used in the present paper. See \cite{DZ98} for a comprehensive treatment of this theory. 

Let $ \Xcal $ denote a topological vector space.  A lower semi-continuous function $ I\colon \Xcal\to [0,\infty] $ is called a {\it rate function\/} if  $ I $ is not identical $ \infty$ and  has compact level sets, i.e., if $ I^{-1}([0,c])=\{x\in\Xcal\colon I(x)\le c\} $ is compact for any $ c\ge 0 $. A sequence $(X_N)_{N\in\N}$ of $\Xcal$-valued random variables $X_N$  satisfies the {\it large-deviation upper bound\/} with {\it speed\/} $a_N$ and rate function $I$ if, for any closed subset $F$ of $\Xcal$,
\begin{equation}\label{LDPupper}
\limsup_{N\to\infty}\frac 1{a_N}\log \P(X_N\in F)\leq -\inf_{x\in F}I(x),
\end{equation}
and it satisfies the {\it large-deviation lower bound\/} if, for any open subset $G$ of $\Xcal$,
\begin{equation}\label{LDPlower}
\liminf_{N\to\infty}\frac 1{a_N}\log \P(X_N\in G)\leq -\inf_{x\in G}I(x).
\end{equation}
If both, upper and lower bound, are satisfied, one says that $(X_N)_N$ satisfies a {\it large-deviation principle}. The principle is called {\it weak\/} if the upper bound in \eqref{LDPupper} holds only for {\it compact\/} sets $F$. A weak principle can be strengthened to a full one by showing that the sequence of distributions of $X_N$ is {\it exponentially tight}, i.e., if for any $L>0$ there is a compact subset $K_L$ of $\Xcal$ such that $\P(X_N\in K_L^{\rm c})\leq {\rm e}^{-LN}$ for any $N\in\N$.

All the above is usually stated for probability measures $\P$ only, but the notion easily extends to {\it sub}-probability measures $\P=\P_N$ depending on $N$. Indeed, first observe that the situation is not changed if $\P$ depends on $N$, since a large deviation principle depends only on distributions. Furthermore, the connection between probability distributions $\widetilde \P_N$ and sub-probability measures $\P_N$ is provided by the transformed measure $\widetilde \P_N(X_N\in A)=\P_N(X_N\in A)/\P_N(X_N\in\Xcal)$: if the measures $\P_N\circ X_N^{-1}$ satisfy a large deviation principle with rate function $I$, then the probability measures $\widetilde \P_N\circ X_N^{-1}$ satisfy a large deviation principle with rate function $I-\inf I$.

One standard situation in which a large deviation principle holds is the case where $\P$ is a probability measure, and $X_N=\frac 1N(Y_1+\dots+Y_N)$ is the mean of $N$ i.i.d.~$\Xcal$-valued random variables $Y_i$ whose moment generating function $M(F)=\int {\rm e}^{F(Y_1)}\,\d\P$ is finite for all elements $F$ of the topological dual space $\Xcal^*$ of $\Xcal$. In this case, the abstract Cram\'er theorem provides a weak large deviation principle for $(X_N)_{N\in\N}$ with rate function equal to the Legendre-Fenchel transform of $\log M$, i.e., $I(x)=\sup_{F\in \Xcal^*}(F(x)-\log M(F))$. An extension to independent, but not necessarily identically distributed random variables is provided by the abstract G\"artner-Ellis theorem.

In our large deviations results we shall rely on the following conventions. For $X=\Omega $ or $X=\Omega_k, k\in\N $, we conceive $ \Mcal_1(X)$  as a closed convex subset of the
space $\Xcal=\Mcal(X)$ of all finite signed Borel measures on $X$. This is a
topological Hausdorff vector space whose topology is induced by the set
$\Ccal_{\rm b}(X)$ of all continuous bounded functions  $X\to\R$. Then $ \Ccal_{\rm b}(X)$ is  the topological dual of $ \Mcal(X)$ \cite[Lemma~3.2.3]{DS01}. When we speak of a
large deviation principle for $\Mcal_1(X)$-valued random variables, then we
mean a principle on $\Mcal(X)$ with a rate function that is tacitly
extended from $\Mcal_1(X)$ to $\Mcal(X)$ with the value $+\infty$.

\subsection{Bose functions}\label{appendix-Bosefunction}
These functions are defined by 
\begin{equation}\label{defBosefunctionapa}
g_{s}(\alpha)=\frac{1}{\Gamma(s)}\int_0^\infty\frac{t^{s-1}}{{\rm e}^{t+\alpha}-1}\d t=\sum_{k=1}^\infty k^{-s}{\rm e}^{-\alpha k}\quad\mbox{ for all } \alpha>0\mbox{ and all } s>0,
\end{equation}
and also $ \alpha=0 $ and $ s>1 $. In the latter case,
$$
g_s(0)=\sum_{k=1}^\infty k^{-s}=\zeta(s),
$$
which is the zeta function of Riemann. The behaviour of the Bose functions about $ \alpha=0 $ is given by
$$
g_s(\alpha)=\left\{\begin{array}{r@{\;,\;}l}
\Gamma(1-s)\alpha^{s-1}+\sum_{k=0}^\infty\zeta(s-k)\frac{(-\alpha)^k}{k!} & s\not= 1,2,3,\ldots\\[1.5ex]
\frac{(-\alpha)^{s-1}}{(s-1)!}\Big[-\log\alpha +\sum_{m=1}^{s-1}\frac{1}{m}\Big]+\sum_{\heap{k=0}{k\not= s-1}}\zeta(s-k)\frac{(-\alpha)^k}{k!} & s=1,2,3,\ldots                   \end{array}\right..
$$
At $ \alpha=0 $, $ g_s(\alpha) $ diverges for $ s\le 1 $; indeed for all $ s$ there is some kind of singularity at $ \alpha=0 $, such as a branch point. For further details see \cite{Gram25}.
\subsection{Quantum statistics for Bosons}\label{qm-sec}
We review some basic facts on quantum systems of non-interacting Bosons. Fix a box $ \L\subset\R^d $. We consider a system of $ N $ identical Bosons enclosed in $ \L $ in thermal equilibrium at inverse temperature $ \beta>0 $. The total energy of the system is given by  a Hamilton operator 
$$ 
H_{\L}^{\ssup{N}}=\sum_{i=1}^N-\Delta_i,
$$ where $ -\Delta_i $ is the Laplacian with some boundary conditions representing the kinetic energy of particle $ i $. The states of the system belong to the {\it symmetrised} subspace $ \Hcal_{\L}^{\ssup{N,{\rm sym}}} $ of the $ N$-fold tensor product $ \Hcal_{\L}^{\otimes N} $ of the one particle Hilbert space $ \Hcal_{\L} =L^2(\L) $. The projection $ P^{\ssup{\rm sym}}_N $ of $ \Hcal_{\L}^{\otimes N} $ onto  $ \Hcal_{\L}^{\ssup{N,{\rm sym}}} $ reads
$$
P^{\ssup{\rm sym}}_N=\frac{1}{N!}\sum_{\pi\in\Sym_N}U_{\pi},
$$
where $ U_{\pi}\colon\Hcal_{\L}^{\otimes N}\mapsto\Hcal_{\L}^{\otimes N} $ is a unitary representation of the permutation group $ \Sym_N $ on $ \Hcal_{\L}^{\otimes N} $ defined by
$$
U_{\pi}(\psi_1\otimes\cdots\otimes\psi_N)=\psi_{\pi(1)}\otimes\cdots\otimes\psi_{\pi(N)},\quad \psi_k\in\Hcal_{\L},\quad k=1,\ldots,N, \pi\in\Sym_N.
$$
The quantum canonical partition function is
$$
Z_{\L}(\beta,N)=\tr_{\Hcal_{\L}^{\ssup{N,{\rm sym}}}}\Big({\rm e}^{-\beta H_{\L}^{\ssup{N}}}\Big)=\tr_{\Hcal_{\L}^{\otimes N}}\Big(P^{\ssup{\rm sym}}_N{\rm e}^{-\beta H_{\L}^{\ssup{N}}}\Big)=\frac{1}{N!}\sum_{\pi\in\Sym_N}f_{\L}(\pi),
$$
where $ f_{\L} $ is a function on $ \Sym_N $ defined by
$$
f_{\L}(\pi)=\tr_{\Hcal_{\L}^{\otimes N}}\Big(U_{\pi}{\rm e}^{-\beta H_{\L}^{\ssup{N}}}\Big),\quad \pi\in\Sym_N.
$$
Each permutation $ \pi\in\Sym_N $ can be written as a product of independent cycles. A cycle of length $ k $ with $ 1\le k\le N $ is a chain of permutations, such as $ 1 $ goes to $2$, $2$ goes to $3$, etc. until $ k-1 $ goes to $ k $ and finally $ k $ goes to $1$. A permutation $ \pi\in\Sym_N $ may be built up of $ r_1 \; 1$-cycles, $r_2\; 2$-cycles,$\ldots, r_N\; N$-cycles; we then say $ \pi $ has the cycle structure or is o type $ (r_1,\ldots,r_N) $ with the 
$ \sum_{k=1}^Nkr_k=N $. Hence, any partition of the integer $ N $ defines a conjugacy class of permutations such that the trace operation for the elements of this class remains constant. The genuine task of quantum statistical mechanics is to prove and analyse the thermodynamic limit $ \lim_{\L\to\R^d}1/|\L|\log Z_{\L}(\beta,N) $ such that $ N/|\L|\to\rho\in (0,\infty) $ (\cite{Rue69}). Feynman 1953 introduced the functional integration methods for traces, see \cite{Gin70} for details. Calculations in the Fourierspace and in the grandcanonical ensemble go back to Bose and Einstein 1925. Bose-Einstein condensation is defined as a macroscopic occupation of the zero mode (in Fourierspace), which corresponds to particles with zero wave vector, i.e, zero momentum. This is only present for systems of Bosons with given density at sufficiently low temperatures. A mathematical criterion for Bose-Einstein condensation has been proposed by Onsager and Penrose \cite{OP56}. Bose-Einstein condensation for non-interacting particles is equivalent to the appearance of long cycles (cycles that grow as the system size) (\cite{S02},\cite{Uel06}). Proofs for Bose-Einstein condensation in the thermodynamic limit have been obtained for mean-field models (\cite{BCMP05},\cite{DMP05}), and for lattice systems with hard-core exclusion interaction in \cite{LSSY05}. For systems confined in trap potentials, Bose-Einstein condensation was proved in the so-called Gross-Pitaevskii limit, which is different to the usual thermodynamic limit (\cite{LSSY05}).


\subsection*{Acknowledgments} The authors thanks Tony Dorlas and Joe Pul\'{e} for discussions on cycles at Dublin Institute for Advanced Studies.

\noindent






\end{document}